\magnification=\magstep1

\font\title=cmr10 scaled 1200

\font\teo=cmcsc10 
 
\font\tenmsam=msam10 
\font\tenmsb=msbm10 
\font\sevenmsb=msbm7 
\font\fivemsb=msbm5 
\newfam\msbfam 
\textfont\msbfam=\tenmsb 
\scriptfont\msbfam=\sevenmsb 
\scriptscriptfont\msbfam=\fivemsb 
\def\quadratino{{\tenmsam\char003}} 
\def\cvd{{\unskip\nobreak\penalty50 
\hskip1em\hbox{ }\hskip2em\hbox{}\nobreak\hfil\quadratino 
\parfillskip=0pt \finalhyphendemerits=0 \par} 
\vskip8pt plus4pt minus2pt} 

\font\tengoth=eufm10 
\font\sevengoth=eufm7 
\font\fivegoth=eufm5 

\newfam\gothfam \scriptscriptfont\gothfam=\fivegoth
\textfont\gothfam=\tengoth \scriptfont\gothfam=\sevengoth
\def\goth{\fam\gothfam\tengoth}

\def\Bbb#1{{\fam\msbfam\relax#1}} 

\def\N{{\Bbb N}} 
\def\Q{{\Bbb Q}} 
\def\Z{{\Bbb Z}} 
\def\C{{\Bbb C}}
\def\R{{\Bbb R}}
\def\Pr{{\Bbb P}}
\def\A{{\Bbb A}}

\def\s{{\cal S}}
\def\c{{\cal C}}
\def\O{{\cal O}}

\def\r{{\cal R}}
\def\F{{\cal F}}

\def\GL{{\rm GL}}
\def\SL{{\rm SL}}

\def\SO{{\rm SO}}
\def\Gm{{\Bbb G}_m}
\def\Ga{{\Bbb G}_a}
\def\T{{\Bbb T}}
\def\Gal{{\cal G}{\rm al}}
\def\PGL{{\rm PGL}}
\def\PSL{{\rm PSL}}
\def\End{{\sl End}}

\def\t{{\goth t}}
\def\z{{\goth z}}
\def\n{{\bf n}}

\def\y{{\goth y}}

\def\u{{\goth u}}
\def\g{{\goth g}}
\def\uno{{\bf 1}_N}

\def\second{{\prime\prime}}
\def\rank{{\rm rank}}
\def\spec{{\rm spec}}

\centerline{{\title Rational fixed points for linear group actions}}

\bigskip

\centerline{\sl Pietro Corvaja}

\bigskip

$${}$$

\bigskip

\noindent{\bf Abstract}. We prove a version of Hilbert Irreducibility Theorem for linear algebraic groups.  Given  a connected linear algebraic group $G$, an affine variety $V$ and a finite map $\pi:V\rightarrow G$, all defined over a finitely generated field $\kappa$ of characteristic zero, Theorem 1.6 provides the natural necessary and sufficient condition under which  the set $\pi(V(\kappa))$ contains a Zariski dense sub-semigroup $\Gamma\subset G(\kappa)$; namely, there must exist an unramified covering $p:\tilde{G}\rightarrow G$ and a map $\theta:\tilde{G}\rightarrow V$ such that $\pi\circ \theta=p$. In the case $\kappa=\Q$, $G=\Ga$ is the additive group, we reobtain the original Hilbert Irreducibility Theorem. 

Our proof uses a new diophantine result, due to Ferretti and Zannier [F-Z]. As a first application, we obtain (Theorem 1.1) a necessary condition for the existence of rational fixed points for all the elements of a Zariski-dense sub-semigroup of a linear group acting morphically on an algebraic variety. A second application concerns the characterisation of algebraic subgroups of $\GL_N$ admitting a Zariski-dense sub-semigroup formed by matrices with at least one rational eigenvalue.
\bigskip

$${}$$

\noindent {\bf \S 1 Introduction}. 

A general principle in the theory of diophantine equations asserts that if an equation admits ``many" rational solutions, there should be a geometric reason explaining such abundance. We consider here a (multiplicative) semigroup of $N\times N$ matrices with rational entries: we suppose that all of them admit rational eigenvalues and deduce the natural geometrical consequences. Such consequences, stated in Theorem 1.2 below, will concern the algebraic group generated by the given semigroup.  Consider the natural action of $\GL_N$ on $N-1$-dimensional projective space $\Pr_{N-1}$: for a non-singular matrix with rational entries, the fact of having a rational eigenvalue amounts to having a rational fixed point in $\Pr_{N-1}$; hence we are naturally led to  consider a linear-group action on an arbitrary algebraic variety. We shall suppose that each element of a given Zariski-dense semigroup has rational fixed points and deduce again the natural geometric consequence (Theorem 1.1).   \smallskip

More precisely, let $\kappa$ be a field of characteristic $0$, finitely generated over the prime field $\Q$. From now on, by {\it rational} we shall mean $\kappa$-rational, unless otherwise stated. Let $X$ be an algebraic variety, and $G$ an algebraic group, both defined over $\kappa$. Suppose that $G$ acts $\kappa$-morphically on $X$ [Bo 2, \S 1.7].     
\smallskip

Our main theorem will be  

\bigskip

\noindent{\bf Theorem 1.1}. {\it Let the finitely generated field $\kappa$, the algebraic group $G$, the variety $X$ and the action of $G$ on $X$ be as above. Suppose moreover that $G$ is connected. Let $\Gamma\subset G(\kappa)$ be a Zariski-dense sub-semigroup. If the following two conditions are satisfied:
\smallskip

\item{(a)} for every element $\gamma\in\Gamma$ there exists a rational point $x_\gamma\in X(\kappa)$ fixed by $\gamma$;

\item{(b)} there exists at least one element $g\in G$ with only finitely many fixed points;
\smallskip

\noindent then  \smallskip

\item{(i)} there exists a rational map $w: G\rightarrow X$, defined over $\kappa$, such that for each element $g$ in its domain, $g(w(g))=w(g)$.

\noindent If moreover $X$ is projective, then 

\item{(ii)} each element $g\in G(\kappa)$ has a rational fixed point in $X(\kappa)$.
}
\medskip

We remark at once that the stronger conclusion that the group $G$ itself admits a fixed point, i.e. the rational map $w$ can be taken to be constant, does not hold in general (see Example 1.8 below). Example 1.8 bis shows that to prove the second conclusion $(ii)$, the hypothesis that the variety $X$ is projective cannot be omitted. On the contrary, we do not know whether hypothesis $(b)$ is really needed for $(i)$ and $(ii)$ to hold.
\medskip

As we mentioned, the starting point of this work was the investigation of semigroups of matrices, with rational entries and rational eigenvalues.
The following result gives a complete classification of such semigroups:
\medskip

\noindent{\bf Theorem 1.2}. {\it Let $\kappa$ be a finitely generated field as above, $1\leq r\leq N$ be two natural numbers. Let $G\subset\GL_N$ be a connected algebraic subgroup defined over $\kappa$,  $\Gamma\subset G(\kappa)$ a Zariski-dense sub-semigroup. Suppose that each matrix $\gamma\in\Gamma$ has at least $r$ rational eigenvalues (counting multiplicities). Then \smallskip

\item{(i)} each matrix $g\in G(\kappa)$ has at least $r$ rational eigenvalues, counting multiplicities; \smallskip

\item{(ii)} there exists an algebraic group homomorphism $G\rightarrow\Gm^r$, defined over $\kappa$, 
$$ 
G\ni g\mapsto (\chi_1(g),\ldots,\chi_r(g))\in\Gm^r 
$$ 
such that for each $g\in G$ the characteristic polynomial of the matrix $g$ is divisible by the degree $r$ polynomial $(T-\chi_1(g))\cdots(T-\chi_r(g))$. } 
\medskip

One of the motivations of the present work arises from a paper of Bernik [Be] concerned with semigroups of matrices whose spectra lie on a fixed finitely generated field. Bernik's result, which for simplicity we state below in a slightly weaker form, will be easily deduced from the case $r=N$ of Theorem 1.2: \medskip

\noindent {\bf Corollary 1.3} (Bernik). {\it Let $\kappa$ be a finitely generated field as before, $\Gamma\subset\GL_N(\kappa)$ be a group of matrices such that each element of $\Gamma$ has its spectrum contained in $\kappa$. Then $\Gamma$ contains a normal solvable subgroup of finite index. } \bigskip

For later convenience, we restate the case $r=1$ of Theorem 1.2 in a slightly stronger form: 

\medskip

\noindent{\bf Corollary 1.4}. {\it Let $\Gamma\subset \GL_N(\kappa)$ be a semi-group of matrices with rational entries and let $G$ be the Zariski-closure of $\Gamma$ (which is an algebraic subgroup of $\GL_N$). Suppose $G$ is connected. The following conditions are equivalent:\smallskip

\item{(i)} every matrix in the semigroup $\Gamma$ has at least one rational eigenvalue;\smallskip

\item{(ii)} every matrix in the group $G(\kappa)$ has at least one rational eigenvalue; \smallskip

\item{(iii)} there exists a character $\chi:G \rightarrow \Gm$, defined over $\kappa$, such that for every $g\in G$, $\chi(g)$ is an eigenvalue of $g$;\smallskip

\item{(iv)} there exists a rational map $w: G\rightarrow\Pr_{N-1}$, defined over $\kappa$, such that for all $g$ in its domain the point $w(g)$ is fixed by the projective automorphism defined by $g$.}

\medskip

Here and in the sequel, by a {\it character} of an algebraic group $G$ we mean an algebraic group homomorphism (see [Bo2, Chap. II, \S 5]), which might be the trivial one. In the case $G$ is semisimple and defined over the reals and $\Gamma\subset G(\R\cap\kappa)$, Corollary 1.4 is contained in a Theorem by Prasad and Rapinchuk  [P-R 2, Theorem 1] (see also [P-R 1, Theorem 2]). Both these works and the one by Bernik [Be] use $p$-adic methods.

As an immediate application of the above corollary, we obtain that each Zariski-dense subgroup of $\GL_N(\kappa)$ or of $\SL_N(\kappa)$ (for $N \geq 2$) contains a matrix with no rational eigenvalue.

\bigskip

An interesting case of Theorem 1.1 arises from the natural action of $\GL_N$ on Grassmannians.  We denote by $\F(r;N)$ the variety of $r$-dimensional subspaces of a fixed $N$-dimensional vector space (say the group variety $\Ga^N$); alternatively, $\F(r;N)$ is the variety of $r-1$-dimensional spaces in $\Pr_{N-1}$. Every algebraic group $G\subset\GL_N$ acts naturally on $\F(r,N)$; an element $g\in G$ fixes a point $\omega\in\F(r;N)$ whenever the subspace $\omega$ is invariant for $g$. More generally, one can consider {flag  varieties}: given integers $0<r_1<\ldots< r_h<N$ let $\F(r_1,\ldots,r_h;N)$ be the variety classifying filtrations $V_1\subset\ldots V_h\subset \Ga^N$ formed by $r_i$-dimensional subspaces $V_i$; in particular, whenever $h=1$ we reobtain the Grassmannian. Of course the group $G\subset\GL_N$ acts naturally on the variety $\F(r_1,\ldots,r_h;N)$ and the fixed points for an element $g\in G$ are just the filtrations of invariant subspaces.
Bernik's theorem (Corollary 1.3) is concerned with the action of an algebraic subgroup $G\subset\GL_N$ on the maximal flag variety $\F(1,2,\ldots,N-1;N)$. 
As a corollary of Theorem 1.1 we obtain the following general statement:
\medskip

\noindent{\bf Theorem 1.5}. {\it Let $0<r_1<\ldots<r_h<N$ be integers as before; let $G\subset\GL_N$ be a connected algebraic group defined over the finitely generated field $\kappa$ as before. Suppose that there exists a matrix $g\in G$ with $N$ distinct eigenvalues. Let $\Gamma\subset G(\kappa)$ be a Zariski-dense semigroup. Suppose that each matrix $\gamma\in\Gamma$ admits a filtration $\{0\}\subset V_1\subset\ldots\subset V_h\subset \kappa^N$ of invariant  subspaces, defined over $\kappa$, with $\dim(V_i)=r_i$. Then:
\smallskip

\item{(i)} every matrix $g\in G(\kappa)$ admits such a decomposition;
\smallskip

\item{(ii)} there exists a rational map $w:G\rightarrow \F(r_1,\ldots,r_h;N)$, defined over $\kappa$, such that for each matrix $g$ in its domain, $w(g)$ is an invariant filtration for $g$. 

\noindent In the case $h=N-1$, so $(r_1,\ldots,r_h)=(1,\ldots,N-1)$, such a map can be taken to be constant.}
\medskip

The condition that some matrix in $G$ has distinct eigenvalues is probably not necessary; for instance, it is not necessary in the case of the complete flag variety ($(r_1,\ldots,r_h)=(1,\ldots,N-1)$).

As a Corollary of Theorem 1.5, we obtain that every Zariski-dense semigroup   $\Gamma\subset\SL_N(\kappa)$ (or $\Gamma\subset\GL_N(\kappa)$) contains a matrix whose characteristic polynomial is irreducible. As an application of the Theorem 1.6 below, we could also prove that its Galois group (over $\kappa$) is infinitely often the full simmetric group on $N$ letters (see Corollary 1.11 for a general statement).

\bigskip

The proof of Theorem 1.2 and its corollaries, including Bernik's theorem, are reduced to certain diophantine equations involving linear recurrences, to be solved in finitely generated groups. Such equations could be dealt with by rather elementary methods, involving height considerations. 

On the contrary, the proof of Theorem 1.1, or even just Theorem 1.5 above, makes use of   completely different techniques from the theory of diophantine equations involving power sums. Such techniques, introduced by Zannier in [Z1] and developed by Ferretti and Zannier in [F-Z], lead to the results stated in \S 3. We choose, for shortness, to use these diophantine results to derive all our main theorems, including Theorem 1.2. 
As a step in our proof,  we shall also obtain the following  
\medskip

\noindent{\bf Theorem 1.6}. {\it Let the field $\kappa$ and the connected linear algebraic group $G$ be as before. Let $V$ be a smooth affine algebraic variety of the same dimension as $G$,  $\pi:V\rightarrow G$ a finite map, both defined over $\kappa$. Let $\Gamma\subset G(\kappa)$ be a Zariski-dense semigroup. If $\Gamma$ is contained in the set $\pi(V(\kappa))$, then there exists an irreducible component $V^\prime$ of $V$ such that the restriction $\pi_{|V^\prime}:V^\prime \rightarrow G$ is an unramified cover. In particular, $V^\prime$ has the structure of an algebraic group over $\kappa$.
}
\medskip

The condition that $V$ is smooth could be avoided, up to rephrasing the conclusion, which would state the existence of an unramified covering $p: \tilde{G}\rightarrow G$ and a morphism $\theta:\tilde{G}\rightarrow V$ with $\pi\circ\theta=p$.

In the case where $G=\Ga$, $\kappa=\Q$, and $\Gamma=\N$ is the semigroup of natural numbers, the above statement is equivalent to the original form of Hilbert Irreducibility Theorem [H, Theorem I, p. 107] (see also [Se, chap. 9, \S 6]): since every connected unramified cover of $\Ga$ has degree one, the conclusion in this case is the existence of a section for $\pi$. Hence our Theorem 1.6 can be viewed as a natural generalization of Hilbert Irreducibility Theorem to linear algebraic groups.

When $\Gamma$ is a cyclic group, or semigroup,  (so in particular $G$ is commutative) our Theorem 1.6 is implicit in the main theorem of [F-Z]. 

Theorem 1.6 is also linked with a conjecture of Zannier [Z2, last page], asking for the same conclusion whenever $G=\Gm^N$, but under the much weaker hypothesis that $\Gamma$ is any Zariski-dense {\it set} of $S$-integral points, for a suitable finite set $S$ of places of $k$. (This last condition, that the elements of $\Gamma$ are $S$-integer points, can be replaced by the assumption that $\Gamma$ is contained in a finitely generated subgroup of $G(\kappa)$). In the one-dimensional case, such problem can be solved using Siegel's theorem on integral points on curves; see [D] or [Z2, Ex III.10].  We note, however, that such a strengthening is not possible for non-commutative groups, as the following example proves: let $V\subset\A^1\times\SL_2$ be the subvariety $\{(y,g)\, :\, y^2={\rm Tr}(g)\}$ (where ${\rm Tr}(g)$ is the trace of the matrix $g$) endowed with the projection $\pi:V\rightarrow\SL_2$. It is immediate to check that the {\it set} $\Gamma:=\pi(V(\Z))$ (i.e. the set of matrices with integral entries, whose trace is a perfect square) is dense in $\SL_2$, and nevertheless $V\rightarrow\SL_2$ is a (connected) {\it ramified} cover of $SL_2$. 
\smallskip

For another result of the type of Theorem 1.6, see also [D-Z, Theorem 1]; actually the techniques introduced in [D-Z] indirectly play a role in the present paper.
\medskip

As we said, our results are connected with (and generalize) Hilbert Irreducibility Theorem, although they do not seem to be a direct consequence of it. To further explore this connection we  need a definition, drawn from [Se, chap. 9, \S 1].
\medskip

\noindent {\bf Definition}. Let $\kappa$ be a field of characteristic zero, $X$ be an irreducible algebraic variety defined over $\kappa$. We say that a subset $\Gamma\subset X(\kappa)$ is $\kappa$-thin if there exists an algebraic variety $Y$ and a morphism $\pi: Y\rightarrow X$ defined over $\kappa$ such that 

\item{$(a)$} $\Gamma$ is contained in $\pi(Y(\kappa))$,

\item{$(b)$} the generic fiber of $\pi$ is finite and $\pi$ has no rational section over $\kappa$.
\medskip

The attentive reader would  note that Serre's definition actually is the above one with $X=\Pr_N$ or $X={\Bbb A}^N$. A generalization of Hilbert Irreducibility Theorem [Se, chap. 9] then states that $\Pr_N(\kappa)$ (and ${\Bbb A}^N(\kappa)$) is not $\kappa$-thin, when $\kappa$ is finitely generated over $\Q$.
\smallskip

Consider for instance Corollary 1.4: let us show a deduction of a weak form of property $(iii)$ in the Corollary from condition $(ii)$, using Hilbert Irreducibility Theorem. 
Let $V\subset G\times \Gm$ be the closed subvariety formed by the pairs $(g,\lambda)$ where $\lambda$ is an eigenvalue of $g$ and let $\pi:V\rightarrow G$ be the projection on the first factor. Note that $V$ has the same dimension as $G$. Condition $(ii)$ is equivalent to the morphism $p$ being surjective on the set of rational points, i.e. to the equality $\pi(V(\kappa))=G(\kappa)$. If $G$ is a $\kappa$-rational variety (it is always rational over a suitable finite extension of $\kappa$, see [Bo2]), the set $G(\kappa)$ is not thin in the above sense; hence the variety $Y$ must have a component on which the morphism $\pi$ has degree one, so $\pi$ admits a rational section $g\mapsto (g,\chi(g))$ defined over $\kappa$: this is conclusion $(iii)$, up to the fact that $\chi$ is a character; this last property, however, follows from general algebraic group theory (see Lemma 4.9).

The main novelity in Corollary 1.4 is that it suffices that $\pi(V(\kappa))$ covers the rational points in a Zariski-dense semigroup to deduce the existence of such a section.
To obtain such a conclusion, one can apply Theorem 1.1 (or Theorem 1.6, as we shall do), although, as we already said, in this particular case the proof is technically simpler than the proof of Theorem 1.1 in  general.

\medskip  

A corollary of Theorem 1.6 reads:
\medskip

\noindent {\bf Corollary 1.7}. {\it Let $G$ be a connected simply connected (linear) algebraic group defined over the finitely generated field $\kappa$. Then no Zariski-dense semigroup $\Gamma\subset G(\kappa)$ is $\kappa$-thin.
}
\medskip

On the other hand, it is clear that non-simply connected linear algebraic groups always admit thin Zariski-dense subgroups. In fact, if $G$ is not simply connected, it admits an unramified cover $\pi:G^\prime\rightarrow G$, with $\deg(\pi)>1$, which is also a morphism of algebraic groups over $\kappa$. Since the group $G^\prime(\kappa)$ is Zariski-dense in $G^\prime$ [Bo2, 18.3], the subgroup $\Gamma:=\pi(G^\prime(\kappa))$ is $\kappa$-thin and Zariski-dense in $G$.
\medskip

As we remarked, every subgroup of $\GL_N$ acts naturally on the projective space $\Pr_{N-1}$. Corollary 1.4 states that if every element of the semigroup $\Gamma\subset G(\kappa)$ has a rational fixed point (where $G$ is connected and $\Gamma$ is Zariski-dense in $G$), then every element of $g\in G(\kappa)$ has a rational fixed point, and such a fixed point can be chosen as the image at $g$ of a given rational map.  The following example shows that such a rational map cannot be always taken to be constant, i.e. there might exist no fixed point in $\Pr_{N-1}$ for the whole group $G$, so no vector which is an eigenvector for each $g\in G$.

\medskip

\noindent{\bf Example 1.8}. Let $\c\subset\Pr_2=X$ be a smooth conic defined over $\kappa$ and let $G$ be the group of projective automorphisms of $\Pr_2$ leaving $\c$ invariant (it is always isomorphic over the algebraic closure $\bar{\kappa}$ of $\kappa$ to the linear group $\PGL_2$). Put $\Gamma=G(\kappa)$ and let $\gamma\in\Gamma$ be one of its elements; we want to prove that $\gamma$ has at least one rational fixed point in $\Pr_2$. Recall that $\c$ is isomorphic over $\bar{\kappa}$ to $\Pr_1$, and each element of $G(\bar{\kappa})$, apart from the identity, has exactly one or two fixed points on $\c(\bar{\kappa})$. If the automorphism $\gamma$, which is defined over $\kappa$, has just one fixed point on $\c(\bar{\kappa})$, such a point must be rational; if otherwise it has two fixed points $x_1,x_2\in\c(\bar{\kappa})$, these points are either both rational or quadratic conjugates. In any case, the pair of tangent lines to $\c$ through $x_1,x_2$ is a singular conic  defined over $\kappa$, and invariant under $\gamma$. Its only singular point, which is the intersection of the two tangent lines, is necessarely  rational, and  is fixed by $\gamma$. So, in any case, $\gamma$ has a rational fixed point, but of course there is no point on $\Pr_2$ which is fixed by the whole action of $G$.
\medskip

A modification of the above example will be used to show that the hypothesis of the projectivity of $X$ cannot be omitted from Theorem 1.1, if one wants the conclusion $(ii)$:
\medskip

\noindent{\bf Example 1.8 bis}. 
 Let $\Gamma\subset\PSL_2(\R)$ be a finitely generated, discrete, Zariski-dense subgroup with no parabolic elements. (Such subgroups exist: for instance they arise as fundamental groups of hyperbolic compact Riemann surfaces. Indeed, such a Riemann surface  can be realized as the quotient of the Poincar\'e upper half-plane by the canonical action of a Fuchsian subgroup $\Gamma$ of $\PSL_2(\R)$. The compacteness of the quotient implies that $\Gamma$ has no parabolic elements.) 
We let $\Gamma$ act both on $\Pr_1$ and on the symmetric square of $\Pr_1$, which is isomorphic to $\Pr_2$, so we embed it into $\PGL_3$.
Its action on $\Pr_2$ preserves the image $\c\subset\Pr_2$ of the diagonal of $(\Pr_1)^2$, which is a smooth conic in $\Pr_2$; we are then in a particular case of the situation of Example 1.8 above. Since $\Gamma$ is a finitely generated group, there exists a finitely generated field $\kappa$ over which every element of $\Gamma$, viewed as an element of $\PGL_3$ (or of $\PSL_2$), is defined. As observed in Example 1.8, each element of $\PSL_2(\kappa)$, hence every element of  $\Gamma$, has a $\kappa$-rational fixed point in $\Pr_2$.  If $\gamma\in\PSL_2(\kappa)$ is not parabolic, then it has a fixed point in the complement of $\c$. Put $X:=\Pr_2\setminus\c$. Since $\Gamma$ has only non-parabolic elements, each $\gamma\in\Gamma$ has a rational fixed point in $X(\kappa)$; neverthless there are (parabolic) elements in $\PSL_2(\kappa)$ which have no (rational) fixed points in $X$ (but only in $\c(\kappa)$).

\medskip

We now show that the hypothesis on the connectedness of the Zariski-closure $G$ of $\Gamma$ in Theorem 1.1 (or in Corollary 1.3) cannot be removed: \medskip

\noindent{\bf Example 1.9}. Let $\kappa=\Q$ be the field of rational numbers. Consider the disconnected subgroup $G$ of $\GL_2$, defined over $\Q$, formed by the matrices of the form $\left(\matrix{\lambda&0\cr 0&\mu}\right)$ or $\left(\matrix{0&\lambda \cr \mu&0}\right)$, for nonzero $\lambda,\mu$. Let $\Gamma\subset G(\Q)$ be the subgroup formed by the matrices of $G(\Q)$ whose entries are squares in $\Q$. It is easily checked that $\Gamma$ is indeed a subgroup, that it is Zariski-dense in $G$ and that each matrix $g\in \Gamma$ has its eigenvalues in $\Q$. Neverthless, $G$ admits no character $\chi$ such that $\chi(g)$ is always an eigenvalue of $g$. \medskip

An application of Theorem 1.5 concerns semigroups of endomorphisms of compact tori $\R^N/\Z^N$, where as usual $\Z^N$ denotes the lattice of integral points in $N$-dimensional Euclidean space $\R^N$. An endomorphism of a torus $\R^N/\Z^N$, viewed as a topological group, is induced by the multiplication on $\R^N$ by a $N\times N$ matrix $g$ with integral entries. Such an endomorphism is surjective if and only if $g\in\GL_N(\Q)$, and is an automorphism if and only if $g$ belongs to $\GL_N(\Z)$. We shall identify $N\times N$ matrices with integral coefficients with endomorphisms of the torus, and shall denote by $\End(\R^N/\Z^N)$ the semigroup of such matrices (or endomorphisms). Let $g\in\GL_N(\Q)$ be such a matrix (so that $g$ has integral coefficients); the invariant circles on the torus $\R^N/\Z^N$ for $g$ correspond to $\Q$-rational fixed points for the projective automorphism that $g$ defines on $\Pr_{N-1}$ (recall that $\Pr_{N-1}$ is the set of lines in a $N$-dimensional space, and that its $\Q$-rational points correspond to lines defined over $\Q$). More generally, $r$-dimensional invariant subtori correspond to rational points in the Grassmannian $\F(r;N)$. 
A particular case of Theorem 1.5 can be restated as follows: \medskip

\noindent {\bf Corollary 1.10}. 
{\it Let $\Gamma\subset\End(\R^N/\Z^N)$ be a semigroup of surjective endomorphisms of the $N$-dimensional torus. Viewing $\Gamma$ as a sub-semigroup of $\GL_N$, let $G$ be its Zariski-closure, and suppose it is connected. Suppose also that at least one matrix $\gamma\in\Gamma$ has distinct eigenvalues. Let $r$ be an integer, with $0<r<N$. The following conditions are equivalent: \smallskip

\item{(i)} each element of $\Gamma$ admits an invariant $r$-dimensional subtorus; \smallskip

\item{(ii)} there exist a rational map $w: G\rightarrow \F(r;N)$, defined over $\Q$, and a character $\chi:G\rightarrow \Gm$, also defined over $\Q$, such that for each point $\gamma\in G(\Q)$ on which $w$ is well defined, $w(\gamma)\in\F(r;N)(\Q)$ is an invariant $r$-dimensional torus; the determinant of the restriction of $\gamma$ to $w(\gamma)$ is $\chi(\gamma)$.} \medskip

Once the semigroup $\Gamma\subset\GL_N(\Q)$ is given, for instance via a finite set of generators, it is easy to check whether conclusion $(ii)$ of Corollary 1.10 holds. Also one can, for each pair $(r,N)$, classify the algebraic groups $G\subset\GL_N$ acting irreducibly on $\Q^N$, satisfying conclusion $(ii)$. For instance, for the pairs $(1,3)$ and $(2,3)$, the corresponding algebraic groups can be reconducted to Example 1.8 (see [B-O]).
\bigskip

A classical application of Hilbert Irreducibility Theorem concerns Galois groups (over number fields) attached to fibers of coverings of algebraic varieties (see [Se, chap. 10]).
Our Theorem 1.6, which, as we said, is a version of Hilbert Irreducibility Theorem, also admits such applications, as showed by the statement below. Let us first introduce a definition: with an algebraic subgroup $G\subset\GL_N$ we associate its {\it characteristic polynomial}, i.e. the characteristic polynomial of its generic element. It is the polynomial $P(T,g)\in\kappa[G][T]$, with coefficients in the ring of regular functions of $G$, given by
$$
P(T,g)=\det(g-T\cdot \uno)=(-1)^N(T^N-{\rm Tr}(g)T^{N-1}+\ldots\pm \det(g)).
$$
In the above formula $\uno$ stands for the unit matrix in $\GL_N$ and the coefficients of the polynomial on the right-hand side are the invariants of the matrix $g$; they are expressed by   regular functions on the group variety $G$.

If, for instance, $G=\SL_N$, it is easily checked that the characteristic polynomial is irreducible and its Galois group  is the full simmetric group on $N$ elements; on the other hand,  the characteristic polynomial of $\SO(3)$, or any group conjugate to it over $\GL_3(\C)$, is reducible, being divisible by $(T-1)$: this fact, remarked also in [B-O], gives a further explanation for Example 1.8.
\medskip

\noindent{\bf Corollary 1.11}. {\it Let $G\subset\GL_N$  be a connected algebraic subgroup, defined over the finitely generated field $\kappa$. Let $P(T,g)\in\kappa[G][T]$ be the characteristic polynomial associated to the algebraic group $G$ and let ${\cal G}$ be its Galois group over $\kappa[G]$. For every Zariski-dense semigroup $\Gamma\subset G(\kappa)$, there exists a matrix $\gamma\in\Gamma$ whose splitting field over $\kappa$ has a Galois group isomorphic to ${\cal G}$. 
}

\bigskip

The paper is organized as follows: The proof of our main theorems will be given in \S 5. The next section, on specializations of finitely generated rings, is purely technical. The key result of that section is a theorem of Masser, enabling to reduce the general case of a finitely generated field to the number field case. In \S 3, we present a new result on exponential diophantine equations which is the main arithmetic ingredient in the proof of Theorem 1.1. It consists of a (sightly more general) reformulation of a recent theorem of Ferretti and Zannier  [F-Z], which in turn generalizes Zannier's solution to the Pisot's $d$-th root conjecture [Z1]. It is the key step for the proof of Theorem 1.6, which in turn will be used to derive all the other statements.
 The geometric tools needed in the proofs of our main theorems will be developped in \S 4. 
 \medskip

{\it Acknowledgments}: The author is pleased to thank Janez Bernik, Andrea Maffei, David Masser and Umberto Zannier for helpful (oral or electronic) conversations. He is also grateful to the referee for pointing out a gap in a previous proof of Theorem 1.2 and for suggesting several improvements on the presentation of the paper.

Part of this work was prepared during the workshop ``Diophantine approximation and heights", held at Erwin Schr\"odinger Institute in Vienna, February-May 2006. The author would like to thank the organizers for the invitation and the ESI for financial support.

\bigskip

\noindent {\bf \S 2. Specializations.}

 Some of our proofs will need results from height theory, which will be used after specializing to number fields.

First of all let us recall some standard notation on absolute values and heights (see also [Se, chap. 2]). Let $L$ be a number field, $M_L$ its set of places  and $S\subset M_L$ a finite subset  containing all the archimedean ones. We choose normalizations, denoted $|\cdot|_\nu$, of the absolute values at every place $\nu\in M_L$ in such a way that the product formula 
$$ 
\prod_{\nu\in M_L} |x |_\nu =1 \eqno(2.1)
$$
 holds for every nonzero $x\in L$ (here the product runs over all the places of $L$) and the absolute logarithmic Weil height reads 
$$ 
h(x)=\sum_{\nu\in M_L}\log(\max\{1,|x|_\nu\}), \eqno(2.2)
$$ 
the sum running over all places of $L$.
We shall say that such absolute values are normalized with respect to $L$. We denote by $\O_S$ the ring of $S$-integers of $L$, i.e. 
$$ 
\O_S:=\{x\in L\, :\, |x|_v\leq 1\, {\rm \ for\ all\ places}\ v\not\in S\}  
$$ 
and $\O_S^*$ the group of $S$-units, i.e. the unit group of $\O_S$. \smallskip
 
\medskip

As in the previous section, $\kappa\subset\C$ denotes a subfield of the field of complex numbers, finitely generated over the field $\Q$ of rational numbers. If it is an algebraic extension of $\Q$, then it is a number field; otherwise, it is  a transcendental extension of some number field $L$, regular over $L$, i.e. a function field over $L$. Let $R\subset \kappa$ be an integrally closed finitely generated subring of $\kappa$; then the subring of $R$ formed by the algebraic numbers in $R$ is a ring of $S$-integers $\O_S\subset L$. The ring $R$ is the ring of regular functions on a (integral model of a) normal irreducible algebraic variety $X$ defined over $\O_S$: $R=\O_S[X]$. Our first lemma is Corollary 7.5, p. 43 in Lang's book [L]. 
\medskip

\noindent{\bf Lemma 2.1}. {\it Let $R$ be an integrally closed finitely generated ring as before. The group of units $R^*$ is finitely generated.} 
\cvd

\smallskip

\noindent Actually the condition that $R$ be integrally closed can be omitted, as it is in [L], the general case following easily from the particular case of integrally closed rings.  
 
\medskip

Every point $x\in X(\bar{L})$ gives rise to a specialization $R\rightarrow \bar{L}$, by putting $f\mapsto f(x)$. Its image is a number field, written $L(x)$, containing $L$. 
For this reason, an algebraic point $x\in X(\bar{L})$ will also be called a {\it specialization}.\medskip

\noindent {\bf Definition 2.2}. 
We shall say that a specialization $x\in X(\bar{L})$ is {\it good} if it is injective on the group  $R^*$. 
\medskip

Our aim is to prove the existence of a ``large" set of good specializations. For this purpose we use a result of Masser [M].  

We imbedd the affine variety $X$ in an affine space $\A^N$, in such a way that its projection on the first $s$ coordinates is a finite map (of course $s=\dim(X)$). 
The imbedding $X\hookrightarrow \A^N$ defines a {\it (logarithmic) height} function enabling to define the height of  every algebraic point in $X(\bar{L})$: see [Se, chap. 2]. Such a height will depend on the given imbedding, but its fundamental property, i.e. the finiteness of points of bounded degree and height, will always hold (Northcott's Theorem [Se, \S 2.4]). We shall speak of a height function $h:X(\bar{L})\rightarrow \R$ to mean the height corresponding to some fixed imbedding as above.
\smallskip
 
Following the notation of  Masser [M], we denote by  ${\cal E}(d,h)$ the set of points $x\in X({\bar L})$ with $[L(x):L]\leq d$ and $h(x)\leq h$ which are NOT good specializations. Also, for every finite set $T\subset X$ we denote by $\omega(T)$ the degree of a minimal hypersurface of $\A^N$, not containing $X$, but containing $T$. The Theorem in [M, \S 5] reads:

\medskip

\noindent{\bf Lemma 2.3}. {\it For each $d$ there exists a number $C$, depending on $R$ and $d$, such that  for every $h\geq 1$,
$$
\omega({\cal E}(d,h))\leq C h^{\rank(R^*)^2}.
$$
}
\cvd
\medskip

We recall that  in the particular case where $X$ is a curve, a much stronger result is known: Bombieri, Masser and Zannier proved a bound for the set of bad specializations, outside a proper Zariski closed set [B-M-Z]. 
From the above Lemma 2.3 we then obtain in particular the existence of infinitely many good specializations, as proved also by Rumely [R]. In this work, however, we shall need the full strength of Masser's Theorem, which goes beyond the existence of  one (or infinitely many) good specializations. We shall use Lemma 2.3 through the following
\medskip

\noindent{\bf Corollary 2.4}. {\it Consider as before the finite map $p:X\rightarrow\A^s$. 
The set $T\subset\A^s(L)$, formed by the points $\alpha\in\A^s$ such that $p^{-1}(\alpha)$ contains at least one  good specialization, is not $L$-thin.}

  \smallskip

\noindent{\it Proof}. We use a counting argument, combining Masser's Theorem with a quantitative version of Hilbert Irreducibility Theorem given in [Se, \S 13.1, Theorem 3] (we warn the reader that the letter $H$ is [Se, \S 13.1] denotes the {\it exponential} height. We shall also make use of it, putting $H(\cdot)=\exp(h(\cdot))$).

\noindent Recall that by Schanuel's Theorem [Se, \S 2.5], for   $H\geq 1$, the number of points in $\A^s(L)$ whose exponential height is $\leq H$ grows asymptotically as $ c_1\cdot H^{(s+1)[L:\Q]}$, for a positive constant $c_1$. On the other hand, the already mentioned quantitative version of Hilbert Irreducibility Theorem (Proposition 1 of [Se, \S 13.1]) asserts that every $L$-thin set contains at most $c_2 H^{(s+1/2)[L:\Q]}\log H$ points of height $\leq H$. 

\noindent Let us denote by $T(H)$ the set of points in $T$ whose height is $\leq H$. 
Hence, to prove that $T$ is not thin it suffices to prove that 

{\it For every positive number $c$ the set $T(H)$ contains at least $cH^{(s+1/2)[L:\Q]}\log H$ points of height $\leq H$, provided $H$ is large enough with respect to $c$.}

\noindent We shall prove this claim.
   
\noindent Each point  $\alpha\in\A^s(L)$ has pre-images via $p$ of degree $\leq d:=\deg(p)$ and height $\leq c_2 H(\alpha)$, for a constant $c_2$ (depending only on the map $p$). Consider the set $p^{-1}(T(H))$; it contains all the points of the form $p^{-1}(a)$, for $a\in \A^s(L)$, $H(a)\leq H$, apart possibly those in a set ${\cal E}(d,\log(c_2H))$. 
Applying Masser's Theorem 2.3, with $d=\deg (p)$, we obtain that the set ${\cal E}(d,\log(c_2H))$  is contained in a hypersurface   ${\cal S}(H)\subset\A^N$, of degree $\leq c_3 \log(H)^{r^2}$ (where $r=\rank(R^*)$), not containing $X$. The projection $p({\cal S}(H))$ on $\A^s$ is still a hypersurface in $\A^s$,  since $p:X\rightarrow \A^s$ is a finite map. Moreover, its degree cannot increase, hence it is still bounded by $c_3 \log(H)^{r^2}$. 
We then obtain that the complement of $T(H)$ (in the set of points of $\A^s(L)$ of height $\leq H$) is contained in a hypersurface of $\A^s$ of degree $\leq c_3(\log(H))^{r^2}$. 
By Bezout's Theorem, the intersection of $p({\cal S}(H))$ with every line in $\A^s$ contains at most $\deg(p({\cal S}(H)))\le c_3(\log(H))^{r^2}$ points. Consider the set of lines parallel to the $x_s$-axis: they are given by the system of equations $x_1=a_1,\ldots,x_{s-1}=a_{s-1}$. For each point $(a_1,\ldots,a_{s-1},a_s)$ of height $\leq H$, the vector $(a_1,\ldots,a_{s-1})$ has also height $\leq H$. Hence the set of points of $\A^s$ of height $\leq H$ is contained in the finite union of lines of equation $x_1=a_1,\ldots,x_{s-1}=a_{s-1}$, for  $(a_1,\ldots,a_{s-1})$ ranging over the set of points in $\A^{s-1}$ of height $\leq H$. By Schanuel's Theorem again, applied to $\A^{s-1}$, there are $\simeq c_1 H^{s[L:\Q]}$ such  lines. Each of them containes at most $c_3(\log(H))^{r^2}$ points outside $T(H)$. Hence the number of points in $T(H)$ is at least $c_4 H^{(s+1)[L:\Q]}-c_1c_3 H^{s[L:\Q]}(\log H)^{r^2}$, which is $>c H^{(s+1/2)[L:\Q]}\log H$ as soon as $H$ is sufficiently large.
\cvd 

\bigskip

\noindent {\bf \S 3. Exponential diophantine equations}.

In the proof of our main theorems, we shall encounter diophantine equations involving linear recurrence sequences. For the general theory of linear recurrence sequences we refer to the survey paper by van der Poorten [vdP] and to the more recent one by Schmidt [Sch]; we just recall here that such sequences are given as function $\y:\N\rightarrow\C$ by an expression of the form 
$$ 
\y(n)=p_1(n)\alpha_1^n+\ldots+p_k(n)\alpha_k^n\eqno(3.1) 
$$
 where the {\it roots} $\alpha_1,\ldots,\alpha_k$ are nonzero pairwise distinct complex numbers and the {\it coefficients} $p_1,\ldots,p_k$ are polynomial functions. Both the roots and the coefficients are uniquely determined by the sequence. The ``{\it degree } in $n$" of $\y$ is the maximum of the degrees of the polynomials $p_i(X)$. We shall sometimes extend the domain of a linear recurrent sequence to the set $\Z$ of all integers, by the same formula.

Due to the explicit expression (3.1), linear recurrence sequences are also named exponential polynomials.

In case $k=1$ and $p_1$ is a nonzero constant we say that $\y$ is a geometric progression; in that case it takes the form $\y(n)=\beta\cdot\alpha^n$, for $\alpha,\beta\in \C^*$. 
\smallskip

The occurrence of linear recurrence sequences in this work is due to the fact that for given matrices $g,h\in\GL_N$, each entry of the product matrix $h\cdot g^n$ is a linear recurrence sequence in $n$. Hence, for instance, the conditions that all matrices of the semigroup generated by two matrices $h,g$ have a rational eigenvalue implies  the existence of rational solutions to certain diophantine equations involving exponential polynomials.
We shall later describe more deeply the relation between algebraic groups and exponential polynomials.

\smallskip

Actually we shall eventually need a generalization of the above notion to  exponential polynomials in several variables. \medskip

As in the previous section, we let $\kappa\subset\C$ be a finitely generated field of characteristic zero. The symbol $\bar{\kappa}$ will denote the algebraic closure of $\kappa$ inside $\C$. 
 We shall consider only exponential polynomial of the form (3.1) with coefficients in $\bar{\kappa}[X]$ and roots in $\bar{\kappa}^*$. The Galois group  $\Gal(\bar{\kappa}/\kappa)$ acts canonically on the ring of such exponential polynomials:  
namely, if $\y$ is an exponential polynomial given by the formula (3.1) and $\sigma\in\Gal(\bar{\kappa}/\kappa)$ then we define $\y^\sigma$ to be the exponential polynomial  
$$
\y^\sigma(n)=\sigma(\y(n))=(p_1^\sigma)(n)(\sigma(\alpha_1))^n +\ldots +(p_k^\sigma)(n)(\sigma(\alpha_k))^n.
$$
Here, for a polynomial $p(X)\in\bar{\kappa}[X]$,  $p^\sigma(X)$ denotes the polynomial obtained by applying the automorphism $\sigma$ to the coefficients of $p(X)$.
\smallskip

We have the following fact:
\smallskip

\noindent{\bf Lemma 3.1}. {\it Let $\kappa\subset\C$ be a field, $U\subset\bar{\kappa}^*$ be a torsion-free multiplicative group which is invariant for Galois conjugation over $\kappa$; let $\y$ be an exponential polynomial of the form (3.1) with roots $\alpha_i\in U$ and coefficients $p_i(X)\in\bar{\kappa}[X]$. The following are equivalent:

\item{(i)} $\y$ is fixed by the Galois group $\Gal(\bar{\kappa}/\kappa)$;

\item{(ii)} the function $\y:\N\rightarrow\C$ takes values in $\kappa$  at each point $n\in\N$;  

\item{(iii)} the function $\y:\N\rightarrow\C$ takes values in $\kappa$ at infinitely many points $n\in\N$.
}

\smallskip

To explain the requirement that the roots  belong to a given torsion-free multiplicative group, recall that the celebrated Skolem-Mahler-Lech Theorem (see for instance [vdP, 3.6.1]) asserts that {\it if a linear recurrence sequence $\y$, with roots in a torsion-free group, has infinitely many zeros, than it vanishes identically}. This fact will be crucial in many parts of the proofs of our statements in this section.
\medskip

\noindent{\it Proof of Lemma 3.1}. The implications $(i)\Rightarrow(ii)\Rightarrow(iii)$ are trivial. Let us prove that $(iii)$ implies $(i)$. Let $\sigma\in\Gal(\bar{\kappa}/\kappa)$ be a Galois automorphism. Since the equation $\y(n)-\y^\sigma(n)=0$ has by hypothesis infinitely many solutions, the left-hand side is identicaly zero by the Skolem-Mahler-Lech theorem. Hence $\y^\sigma=\y$ as wanted.  
\cvd

\medskip
 
\noindent{\bf Definition}. Let $\kappa,U$ be as in the above Lemma and suppose moreover that $U$ is finitely generated. We let $\r_{\kappa,U}$ be the ring of exponential polynomials satisfying the equivalent conditions of Lemma 3.1.
\smallskip

The ring $\r_{\kappa,U}$ turns out to be a domain.
It is  isomorphic to the $\kappa$-algebra $\kappa[\Ga\times\T]$, for a suitable $\kappa$-torus $\T$. Such a torus is split, i.e. isomorphic to $\Gm^r$ (where $r=\rank (U)$), if and only if $U\subset\kappa^*$.   We notice at once that, for a field extension $\kappa^\prime/\kappa$, the tensor product $\r_{\kappa,U}\otimes_\kappa\kappa^\prime$ is isomorphic to the ring $\r_{\kappa^\prime,U}$. In particular 
$$
\r_{\kappa,U}\otimes_\kappa\kappa(U)\simeq\kappa(U)[\Ga\times\Gm^r]=\kappa(U)[X,T_1,\ldots,T_r,T_1^{-1},\ldots,T_r^{-1}].
$$
 
\smallskip

The first  result  in this section (Theorem 3.2) concerns  algebraic equations involving exponential polynomials to be solved in the finitely generated field $\kappa$. Theorem 3.5 is its natural generalization to several variables.

In the sequel, we shall often write $\y(n)$ to denote the function $\y:\N\rightarrow \kappa$ sending $n\mapsto \y(n)$. Accordingly, we shall also write $f(T,n)$ to denote a polynomial in $T$, with coefficients in the ring $\r_{\kappa,U}$: this notation is justified by the fact that the elements of the ring $\r_{\kappa,U}$ can be viewed as $\kappa$-valued functions in the variable $n\in\N$. Of course, whenever $n$ is a given natural number, then $f(T,n)$ will be a polynomial in $\kappa[T]$. We think that this ambiguity creates no problem, since it will be clarified by the context.
\smallskip

The following statement will be derived from  Theorem 1.1 of [F-Z] (which will be explicitely stated later as Proposition 3.6):
\medskip

\noindent{\bf Theorem 3.2}. {\it Let the finitely generated field $\kappa$  and the finitely generated torsion-free multiplicative group $U\subset\bar{\kappa}^*$ be as before. Let
$f(T,n)\in \r_{\kappa,U}[T]$ be a monic polynomial of the form
$$
f(T,n)=T^d+\y_1(n)T^{d-1}+\ldots+\y_d(n)\eqno(3.2)
$$
where $\y_1,\ldots,\y_d\in\r_{\kappa,U}$ are exponential polynomials.

\noindent Suppose that for every positive integer $n\in\N$ there exists a solution $t\in\kappa^*$ to the equation 
$$
f(t,n)=0.\eqno(3.3)
$$ 
Then there exists an exponential polynomial $\t\in\r_{\kappa,U}$ such that for each $n\in\N$ with $n\equiv 0$ (mod $d!$):
$$
f(\t(n),n)=0.
$$
}
\medskip

In general the condition $n\equiv 0$ (mod $d!$) cannot be omitted:
\smallskip

\noindent{\bf Example 3.3}. Let $\kappa=\Q$ be the field of rational numbers, $U=\{2^n\, :\, n\in\Z\}$ be the cyclic multiplicative group generated by the integer $2$.
Take for $f$ the reducible polynomial 
$$
f(T,n):=T^4-3\cdot 2^n T^2+2\cdot 4^n =(T^2-2^n)(T^2-2\cdot 2^n).
$$ 
Clearly, it has a root in $U\cap\kappa^*=U$ for all choice of an integer $n\in\Z$, but no functional root in $\r_{\Q,U}$. On the contrary, it has four ``functional solutions"  $u_1:\ n\mapsto\sqrt{2}^n, \quad u_2:\ n\mapsto -\sqrt{2}^n$, $u_3:\ n\mapsto\sqrt{2}\sqrt{2}^n$ and $u_4:\ n\mapsto -\sqrt{2}\sqrt{2}^n$ in the ring $\r_{\Q(\sqrt{2}),U^\prime}$, where $U^\prime$ is generated by $\sqrt{2}$. If one replaces $n$ by $2n$ in $f(T,n)$, then one obtains functional solutions already in $\r_{\Q,U}$. 
\medskip

As we mentioned, we shall actually need a generalization to linear recurrence sequences in several variables. Neverthless, the main technical points of the proof appear already in the one variable case. From the above Theorem 3.2 we shall deduce quite formally its natural generalization  to several variables (Theorem 3.5 below). In order to state such a generalization we need some more notation. Let $h\geq 1$ be an integer; we denote by $\r_{\kappa,U}^{\otimes h}$ the ring of polynomial exponential functions in $h$ variables, with roots in $U$. It can be formally defined as the $h$-fold  tensor product 
$$
\r_{\kappa,U}^{\otimes h}=\r_{\kappa,U}\otimes_\kappa\ldots\otimes_\kappa \r_{\kappa,U};
$$
its elements are polynomial exponential functions in $h$ variables, i.e. expressions of the form
$$
\y(n_1,\ldots,n_h)=\sum_{j=1}^k p_j(n_1,\ldots,n_h)\alpha_{j,1}^{n_1}\cdots\alpha_{j,h}^{n_h},
$$
for polynomials $p_1,\ldots,p_h\in \bar{\kappa}[X_1,\ldots,X_k]$ and roots $\alpha_{j,i}\in U$. Again, we require the invariance under Galois conjugation over $\kappa$.

According to the previous notation, we shall write $f(T,\n)$ to denote a polynomial in $T$ with coefficients in a ring $\r_{\kappa,U}^{\otimes h}$ of exponential polynomials in $\n=(n_1,\ldots,n_h)\in\N^h$.

We shall repeatedly use the following remark:
\medskip

\noindent{\bf Remark 3.4}. {\it With the above notation,  the ring $\r_{\kappa,U}^{\otimes h}$ is an integral domain, isomorphic to the ring $\kappa[\Ga^h\times\T^h]=\kappa[H]$, for a suitable $\kappa$-torus $\T$, i.e. to the ring of functions of a commutative algebraic group $H$. For a positive integer $D\in\Z$, the map $\Z^h\ni (n_1,\ldots,n_h)\mapsto (D n_1,\ldots,D n_h)\in\Z^h$ induces an isogeny $H\rightarrow H$ defined (in multiplicative notation) by $\gamma\mapsto \gamma^D$. The units of the ring $\r_{\kappa,U}^{\otimes h}$ are of the form $\u(\n)=\beta\alpha_1^{n_1}\cdots\alpha_h^{n_h}$ for $\beta\in\kappa^*$ and $\alpha_1,\ldots,\alpha_h\in U\cap\kappa^*$.
}
\cvd
 
\medskip

The above mentioned generalization of Theorem 3.2 is the following

\medskip

\noindent{\bf Theorem 3.5} {\it Let $\kappa,U$ be as above, $h$ a positive integer, $f(T,\n) \in\r_{\kappa,U}^{\otimes h}[T]$ a monic polynomial with coefficients in the ring of exponential polynomials in $h$ variables.  If for each vector $(n_1,\ldots,n_h)\in\N^h$ the equation
$$
f(t,n_1,\ldots,n_h)=0\eqno(3.5)
$$
has a rational solution $t\in\kappa$, there exists an exponential polynomial $\t\in\r_{\kappa,U}^{\otimes h}$ such that identically
$$
f(\t(n_1,\ldots,n_h),d!\cdot n_1,\ldots, d!\cdot n_h)\equiv 0.
$$
}

\medskip

One could prove, under the hypotheses of Theorem 3.5, the existence of a functional solution to equation (3.5), i.e. an exponential polynomial $\t^\prime$ satisfying $f(\t^\prime(\n),\n)\equiv 0$ identically (without restricting to the subgroup $d!\cdot\Z^h\subset\Z^h$). Neverthless, such an exponential polynomial might not exist in the ring $\r_{\kappa,U}^{\otimes h}$ (see Example 3.3).
\smallskip

We now state, with our notation, a particular case of Theorem 1.1 in [F-Z]:
\medskip

\noindent{\bf Proposition 3.6}. {\it Let $L\subset\C$ be a number field, $U\subset\bar{L}^*$ a finitely generated torsion-free multiplicative group, invariant under Galois conjugation over $L$. Let $f(T,n)\in\r_{L,U}[T]$ be a monic polynomial with coefficients in the ring of exponential polynomials. Suppose that for all but finitely many integers $n\in\N$, the diophantine equation 
$$
f(t,n)=0
$$
has a rational solution $t\in L$. Then there exists an exponential polynomial $\t$ with algebraic roots and algebraic coefficients, such that identically $f(\t(n),n)\equiv 0$.
}
\medskip
\cvd

There are three main differences between the result above and our statements: (1) we claim in Theorems 3.2,  3.5 that the functional solutions belong to the same ring $\r_{\kappa,U}$, up to restricting to the   arithmetic progression $n\equiv 0$ (mod $d!$); (2) we work with exponential polynomials in several variables; (3) we do not suppose that our finitely generated field $\kappa$ is a number field.

We begin by solving the first problem, via the following three elementary lemmas; they will allow to deduce Proposition 3.10 below from Proposition 3.6; Proposition 3.10 will represent the number field case of Theorem 3.5 for $h=1$:
\medskip

\noindent{\bf Lemma 3.7}. {\it Let $\bar{\kappa}\subset\C$ be an algebraically closed field, $H,H^\prime$ two connected commutative algebraic groups defined over $\bar{\kappa}$ of the same dimension. If $p:H^\prime\rightarrow H$ is an isogeny of degree $d$, then $H$ is isomorphic to $H^\prime$ and there exists an isogeny $\rho:H\rightarrow H^\prime$ such that $p\circ\rho$ is the isogeny sending $\gamma\mapsto\gamma^d$.
}
\smallskip

\noindent{\bf Remark}. In algebraic terms,   Lemma 3.7 states that an unramified (\'etale) extension  of the $\bar{\kappa}$-algebra $\bar{\kappa}[X,T_1,T_1^{-1},\ldots,T_r, T_r^{-1}]$ of degree $d$ is always contained in  the extension of the form $\bar{\kappa}[X,T_1,T_1^{-1},\ldots,T_r,T_r^{-1}][T_1^{1/d},\ldots,T_r^{1/d}]$, which has degree $d^r$.\smallskip

\noindent{\it Proof of Lemma 3.7}. 
We know that $H$ and $H^\prime$ are isomorphic to products of powers of the additive group $\Ga$ and the multiplicative group $\Gm$. Since they are isogenous, they must both be of the  form $\Ga^e\times\Gm^r$, for the same exponents $e,r$, so in particular they are isomorphic. Let us first consider the case $e=0$, i.e. $H=H^\prime=\Gm^r$.
Now each isogeny $p:\Gm^r\rightarrow\Gm^r$ is given by an expression of the form
$$
\pi(t_1,\ldots,t_r)=(t_1^{a_{11}}\cdots t_r^{a_{1r}},\ldots,t_1^{a_{r1}}\cdots t_r^{a_{rr}}),
$$
for a non-singular matrix $A=(a_{ij})_{1\leq i,j\leq r}$ with integral entries.
The degree $d$ of the isogeny $p$ is the absolute value of the determinant of the matrix $A$. Then the matrix  $B:=d\cdot A^{-1}=(b_{i,j})_{1\leq i,j\leq r}$ has integral coefficients, hence determines an isogeny $\rho:\Gm^r\rightarrow\Gm^r$ as above.  Clearly the product $p\circ\rho$ is the isogeny raising to the $d$-th power every element of $H=\Gm^r$.

Let us now consider the general case $H=H^\prime=\Ga^e\times\Gm^r$. We shall use coordinates $(x,y)$ for points of $H,H^\prime$, with $x\in\Ga^e$ and $y\in\Gm$. 
Since $\Ga$ is simply connected, every connected unramified cover of the product $\Ga^e\times\Gm^r$ must isomorphic (as a cover) to a product of the trivial cover of $\Ga^e$ by a  connected unramified cover of $\Gm^r$. Hence there exist an automorphism $\phi$ of $H^\prime=\Ga^e\times\Gm^r$ such that $p\circ\phi (x,y)=(x,p^\prime(y))$, where $p^\prime:\Gm^r\rightarrow\Gm^r$ is an isogeny of degree $d=\deg(p)$. By what we have just proved, there exists an isogeny $\rho^\prime:\Gm^r\rightarrow\Gm^r$ such that the product $p^\prime\circ\rho^\prime$ sends $y\mapsto y^d$. Putting $\rho(x,y):=\phi (dx,\rho^\prime(y))$ we obtain that $p\circ\rho$ is the required isogeny.
\cvd
\medskip

\noindent {\bf Lemma 3.8}. {\it Let $U\subset \tilde{U}\subset\bar{\kappa}^*$ be  finitely generated torsion free multiplicative groups; let $\t\in\r_{\bar{\kappa},\tilde{U}}$ be  an exponential polynomial. If $\t$ is algebraic over the (quotient field of the) ring $ \r_{\bar{\kappa},U}$, then the group generated over $U$ by the roots of $\t$ has the same rank as $U$.
}
\smallskip

\noindent{\it Proof}. We decompose $\tilde{U}$ as a direct sum $\tilde{U}=U_1\oplus U_2$, where $U$ has finite index in $U_1$ and $U_2\cap U=\{0\}$.  Then, letting $s$ be the rank of $U_2$, the ring $\r_{\bar{\kappa},\tilde{U}}$ is isomorphic to $\r_{\bar{\kappa},U_1}[T_1,\ldots,T_s,T_1^{-1},\ldots,T_s^{-1}]$; in such an isomorphism, the elements of $\r_{\bar{\kappa},U}$ are sent to elements of $\r_{\bar{\kappa},U_1}$, so
the ring $\r_{\bar{\kappa},U}$ is identified to a subring of $\r_{\bar{\kappa},U_1}$. Clearly every Laurent polynomial in the ring  $\r_{\bar{\kappa},U_1}[T_1,\ldots,T_s,T_1^{-1},\ldots,T_s^{-1}]$ is transcendental over the subring $\r_{\bar{\kappa},U_1}$, unless it is constant with respect to $T_1,\ldots,T_s$. Hence, it $\t\in\r_{\bar{\kappa},\tilde{U}}$ is algebraic over $\r_{\bar{\kappa},U}$ (so {\it a fortiori} over $\r_{\bar{\kappa},U_1}$) it must have its roots in $U_1$.
\cvd

\medskip

\noindent{\bf Lemma 3.9}. {\it Let $\kappa\subset\C$ be any field, $U\subset\bar{\kappa}^*$ be a finitely generated torsion-free multiplicative group, invariant for the Galois action over $\kappa$. Let $f(T,n)\in\r_{\kappa,U}[T]$ be a monic polynomial of degree $d\geq 1$. Suppose there exists an exponential polynomial $\t$, with arbitrary roots and coefficients, such that identically $f(\t(n),n)\equiv 0$. Then there exists an exponential polynomial $\z\in\r_{\bar{\kappa},U}$, with roots in $U$, such that identically $f(\z(n),d!\cdot n)\equiv 0$.
}
\smallskip

\noindent{\it Proof}. Embedding $\kappa$ in its algebraic closure $\bar{\kappa}$ we can view the polynomial $f(T,n)$ as having coefficients in the ring $\r_{\bar{\kappa},U}$; recall that this ring corresponds to the $\bar{\kappa}$-algebra of a split connected commutative algebraic group of the form $H:=\Ga\times\Gm^r$ (with $r={\rm rank}(U)$). Hence we can think of $f(T,n)$ as a polynomial in $\bar{\kappa}[H][T]$, having its coefficients in the ring of regular functions on $H$; we shall write accordingly  $f(T,g)$ to denote such polynomial. Letting $V:=\spec(\bar{\kappa}[H][T]/(f(T,g))$ be the algebraic variety  $V\subset\A^1\times H$ defined by the equation $f(t,g)=0$. It is endowed with its natural projection $\pi:V\rightarrow H$, which is a morphism of degree $d=\deg(f)$.

Suppose there exists a functional solution $\t$ to the equation $f(t,n)\equiv 0$, as in the statement of the lemma; let $U^\prime$ be the multiplicative group generated by its roots and let $D$ be the order of the torsion part of the group generated by the two groups $U,U^\prime$. Then the exponential polynomial $\t^\prime$ defined by $\t^\prime(n):=\t(Dn)$ has its roots in a finitely generated torsion-free group $U^\second\supset U$. The ring $\r_{\bar{\kappa},U^\second}$ corresponds to the $\bar{\kappa}$-algebra of a connected commutative group $H^\prime$ of the form $H^\prime:=\Ga\times\Gm^{r^\prime}$, where $r^\prime={\rm rank}(U^\second)$; also the inclusion $\r_{\bar{\kappa},U}\hookrightarrow\r_{\bar{\kappa},U^\second}$ corresponds to a surjective algebraic group homomorphism $p:H^\prime\rightarrow H$. Since $\t^\prime$ is algebraic over $\r_{\bar{\kappa},U}$, by Lemma 3.8 we can take $U^\second$ such that the index $[U^\second:U]$ is finite. Hence $r^\prime=r$, $H^\prime=H$ and $p:H\rightarrow H$ is an isogeny. Also, the exponential polynomial $\t^\prime\in\r_{\bar{\kappa},U^\second}$ gives a morphism $\theta:H\rightarrow V$ with $p\circ\theta=\pi$.

We want to bound the degree of minimal  isogeny $p:H\rightarrow H$ for which there exists a morphism $\theta:H\rightarrow V$ as above. Let $V^\prime$ be the image $\theta(H)$ of $H$ in $V$; it is an irreducible variety. Then we have the double inclusion of $\bar{\kappa}$-algebras:
  $p^*(\bar{\kappa}[H])\subset \theta^*(\bar{\kappa}[V^\prime])\subset\bar{\kappa}[H]$. 
We identify the ring $p^*(\bar{\kappa}[H])$ to the ring $\bar{\kappa}[X,T_1,T_1^{-1},\ldots,T_r,T_r^{-1}]$.
Let $d^\prime$ be the degree of the extension $[\theta^*(\bar{\kappa}[V^\prime]):  p^*(\bar{\kappa}[H])]$; $d^\prime$ is the degree of an irreducible factor of the polynomial $f(T,n)\in\r_{\kappa,U}[T]$, so in particular   $d^\prime\leq d$. Since the ring $\theta^*(\bar{\kappa}[V^\prime])$ is contained in  $\kappa[H]$, it is also contained, by Lemma 3.7 (see in particular the remark after the lemma), in an extension  of the form $\bar{\kappa}[X,T_1,T_1^{-1},\ldots,T_r,T_r^{-1}][T_1^{1/e},\ldots,T_r^{1/e}]$ for a suitable integer $e$. We claim it is  contained in an extension of that form with $e=d^\prime$. This last fact can be proved by an argument from Galois theory: namely, the field extension 
$ \bar{\kappa}(X,T_1,\ldots,T_r,)(T_1^{1/e},\ldots,T_r^{1/e})/
\bar{\kappa}(X,T_1 ,\ldots,T_r)$ is Galois with Galois group isomorphic to $(\Z/e\Z)^r$. Since the intemediate field $\theta^*(\bar{\kappa}(V^\prime))$ has degree $d^\prime$ over $\bar{\kappa}(X,T_1,\ldots,T_r,)$, $d^\prime$ must divide $e^r$ and the field $\theta^*(\bar{\kappa}(V^\prime))$ must be the fixed field of a subgroup of $(\Z/e\Z)^r$ containing the subgroup of the multiple of $d^\prime$. Hence it is contained in the sub-extension $\bar{\kappa}(X,T_1,\ldots,T_r,)(T_1^{1/d^\prime},\ldots,T_r^{1/d^\prime})$ as claimed. Since $d^\prime$ divides $d!$, it is also contained in the corresponding extension with $e=d!$. This concludes the proof of the lemma. 
\cvd
\medskip

We finally arrive at the number-field version of Theorem 3.5, for $h=1$;
\medskip

\noindent{\bf Proposition 3.10}. {\it Let $L$ be a number field, $U\subset\bar{L}^*$ be a torsion-free multiplicative group, invariant for the action of $\Gal(\bar{L}/L)$. Let $f(T,n)\in\r_{L,U}[T]$ be a monic polynomial of degree $d\geq 1$ with coefficients in the ring $\r_{L,U}$. Suppose that for all $n\in\N$, the equation $f(t,n)=0$ has a rational solution $t\in L$. Then there exists an exponential polynomial $\z\in\r_{L,U}$ such that identically $f(\z(n),d!\cdot n)\equiv 0$.} 

\smallskip

\noindent{\it Proof}. We  first view the polynomial $f(T,n)\in\r_{L,U}[T]\subset\r_{\bar{L},U}[T]$ as having coefficients in $\r_{\bar{L},U}$. Remember that the latter is the $\bar{L}$-algebra $\bar{L}[\Ga\times\Gm^r]$, where $r=\rank (U)$.
By Proposition 3.6 and the above Lemma 3.9 there exists an exponential polynomial $\t\in\r_{\bar{L},U}$, satisfying the equation $f(\t(n),d!\cdot n)\equiv 0$. Its coefficients are polynomials in $\bar{L}[X]$.

The proof will be finished once we know that for at least one solution to $f(\t(n),d!n)\equiv 0$, the exponential polynomial $\t$ lies in $\r_{L,U}$. Suppose not; we wish to obtain a contradiction. First notice that the Galois group $\Gal(\bar{L}/ L)$ acts on the set of such solutions, which is so invariant. For every $\t\in\r_{\bar{L},U}$, let $r(\t)$ be the multiplicity of the solution $\t$ to the equation $f(\t(n),d!n)\equiv 0$.
Then the product
$$
h(T,n):=\prod_{ \t}(T-\t)^{r(\t)} \in\r_{L,U}[T],
$$
ranging over all the solutions $\t$, divides $f(T,d!n)$ in the ring $\r_{L,U}[T]$. (Observe that it has indeed its coefficients in $\r_{L,U}$, not only in $\r_{\bar{L},U}$, because of the invariance of the set of solutions under Galois conjugation).
By Lemma 3.1, each such exponential polynomial takes  values in $L$ at only finitely many integral points $n\in\N$. 
Now, putting 
$$
g(T,n)={f(T,d!n)\over h(T,n)}\in\r_{L,U} [T]
$$
we arrive at a contradiction with Proposition 3.6: the equation $g(t,n)=0$ has, for all large $n$, a rational solution, nevertheless it admits no functional solution, even with arbitrary algebraic roots and coefficients.   \cvd

\medskip

Our next goal is to pass from number fields to arbitrary finitely generated fields.
For this purpose, we shall use the results of the preceding paragraph, especially the specialization Lemma 2.3 and its corollary.
\medskip

Let $U$ be a finitely generated torsion-free multiplicative group of rank $r$. Given a basis $(u_1,\ldots,u_r)$, it can be identified with the group $\Z^r$. For an element $\alpha\in U$, written as  $\alpha=u_1^{a_1}\cdots u_r^{a_r}$, the {\it height} of $\alpha$, with respect to the given basis, is by definition the integer $|a_1|+\ldots+|a_r|$. Of course, for each given basis and each number $K$, there exist only finitely many elements of $U$ of height $\leq K$, with respect to the given basis. Also, for every basis, the height of a product is bounded by the sum of the heights of the factors.

\medskip

\noindent{\bf Lemma 3.11}. {\it Let  $U\subset\C^*$ be a finitely generated torsion-free group, endowed with a basis. Let $\y_1,\ldots,\y_d$ be exponential polynomials with roots in $U$ and consider the monic polynomial
$$
g(T,n)=T^d+\y_1(n)T^{d-1}+\ldots+\y_d(n).
$$
Let $K$ be an integer larger then the height of each root and the degree in $n$ of each exponential polynomial  $\y_1,\ldots,\y_d$.  
Suppose $\t$ is an exponential polynomial, with roots in $U$, satisfying $g(\t,n)\equiv 0$. Then the height of each root of $\t$ is $\leq K$. Also, the degree in $n$ of $\t$ is also bounded by $K$.
}
 
\smallskip

\noindent {\it Proof}. Representing the elements of $U$ as vectors in $\Z^r$, consider the convex hull of the set of roots of $\t$. Let $\alpha$ be a vertex of  this convex set which is also of maximal height. Then the root $\alpha^d$ does appear in the expansion of $\t^d$. Since it simplifies in $g(\t,n)$, it must be equal to some product $\alpha_1\cdots\alpha_j$ of $j\leq d-1$ roots of $\t$ and a root $\beta$ of some $\y_i$. The height of each root $\alpha_1,\ldots\alpha_j$ is bounded by the height of $\alpha$, while $\beta$ has height $\leq K$; then the height of $\alpha$ must also be $\leq K$. Having bounded the height of each vertex by $K$, we obtain the same bound $K$ for the  convex hull, containing all our roots. The argument to bound the degree in $n$ follows the same pattern (but is simpler since  the totally ordered semigroup $\N$ replaces the group $\Z^r$).
\cvd 
\bigskip

We are now ready to prove Theorem 3.5, at least in the particular case where (1)  $h=1$ and (2) $U\subset\kappa^*$.
\medskip

\noindent{\bf Proposition 3.12}. {\it Let $\kappa$ be a finitely generated field, $U\subset{\kappa}^*$ be a torsion-free multiplicative group. Let $f(T,n)\in\r_{\kappa,U}[T]$ be a monic polynomial of degree $d\geq 1$. Suppose that for each $n\in\N$ the equation $f(t,n)=0$ has a solution in $\kappa$. Then there exists an exponential polynomial $\t\in\r_{\kappa,U}$ such that identically $f(\t(n),d!\cdot n)\equiv 0$.
}
\smallskip

\noindent{\it Proof}. Write $f(T,n)=T^d+\y_1(n)T^{d-1}+\ldots+\y_d(n)$; each exponential polynomial $\y_i(n)$ decomposes as a sum $\y_i(n)=p_{1,i}(n)\alpha_{1,i}^n+\ldots+p_{k_i,i}(n)\alpha_{k_i,i}^n$.
 Let $R\subset\kappa$ be a finitely generated integrally closed ring containing all the coefficients of all the polynomials $p_{i,j}$ with $i\in\{1,\ldots,d\}, \, j\in\{1,\ldots,k_i\}$. Suppose also that its group of units $R^*$ contains $U$. Then $R$ is the ring of regular functions on some affine normal variety $X$ over some ring of $S$-integers of a number field $L$ (see \S 2). Accordingly, the exponential polynomials $\y_i(n)$ can also be viewed, for each $n$, as regular functions on $X$; alternatively, the elements of $\r_{\kappa,U}$ can be viewed as functions on $\N\times X$.  For $x\in X(\bar{L})$ we shall consequently write $\y_i(n)(x)$ to denote the corresponding specialized exponential polynomial (or, if $n$ is a given number in $\N$, $\y_i(n)(x)$ will be the corresponding algebraic number).
The exponential polynomial $\y_i(n)(x)$ will have its roots in $U(x):=\{\alpha(x)\, :\, \alpha\in U\}\subset L(x)^*$. (Here, as usual, $L(x)\subset\bar{L}$ denotes the residue field of at the point $x\in X(\bar{L})$).
We also write $f(T,n)(x)$ to denote the polynomial 
$$
f(T,n)(x):=T^d+\y_1(n)(x)T^{d-1}+\ldots+\y_d(n)(x).
$$
Let $x\in X(\bar{L})$ be a given point (i.e. a specialization). Then Lemma 3.10, applied with $L(x)$ instead of $L$, assures that there exists  a functional solution $\t_x\in\r_{L(x),U(x)}$ to the equation $f(t,d!\cdot n)(x)=0$.

Suppose from now on that $x$ is  a {\it good specialization}, i.e. one  that is injective on $R^*\supset U$ (see Def. 2.2). Observe that the basis $u_1,\ldots,u_r$ of $U$ gives by specialization a basis $u_1(x),\ldots,u_r(x)$ of $U(x)$; also, the specialization map $U\rightarrow U(x)$ preserves the height (with respect to these bases). Since the heights of the roots of $\t_x(n)$ are bounded independently of $x$, by the previous lemma, we have only finitely many possibilities for the elements $\alpha\in U$ such that for some good specialization $x$, $\alpha(x)$ is a root of $\t_x$. For the same reason, the degree in $n$ of the polynomial coefficient of $\t_x$ is also bounded independently of $x$. Let $K$ be a bound for the degree in $n$ of $\t_x$ and for the height of its roots, uniform in $x\in X$, and let $\alpha_1,\ldots,\alpha_l$ be all the elements of $U$ having height $\leq K$. 
\smallskip

We search for an exponential polynomial $\t\in\r_{\kappa,U}$ of the form 
$$
\t(n)=p_1(n)\alpha_1^n+\ldots+p_l(n)\alpha_l(n)
$$
where $\deg p_i\leq K$ for $i=1,\ldots,l$, such that the $f(\t(n),d!\cdot n)=0$. We can view the coefficients of $p_1,\ldots,p_l$ as unknowns (in the affine space $\A^{l(K+1)}(\kappa)$) and the condition $f(\t(n),d!\cdot n)\equiv 0$ corresponds to a system of algebraic equations in $\A^{l(K+1)}\times X$. Let $W\subset A^{l(K+1)}\times X$ be the variety associated to such equation. Finally, let $\pi:W\rightarrow X$ be the projection on the second factor; note that $\pi$ has finite degree $\leq d$, since for no specialization $x\in X(\bar{L})$ there can exist more than $d$ functional solutions $\t_x$ to the equation $f(\t_x,d!n)(x)\equiv 0$. The existence of a solution $t_x$ for each good specialization $x\in X(\bar{L})$ means that there exists a point $w\in W(L(x))$ with $\pi(w)=x$.
Our aim (i.e. the proof of the existence of a solution $\t\in \r_{\kappa,U}$) amounts to   proving that $\pi$ admits a section, defined over $L$; in other words, the unknown coefficients of $p_1,\ldots,p_l$ should be given by regular functions on $X$. For this purpose, we shall use Hilbert Irreducibility Theorem.

\noindent Now consider, as in Corollary 2.4, a finite map $p:X\rightarrow \A^s$ ($s=\dim(X)$), defined over $L$.  We shall consider the points of $X(\bar{L})$ which are pre-images of some $L$-rational point in $\A^s$. Due to Corollary 2.4, the set  
$T\subset\A^s(L)$ defined by 
$$
T:=\{\alpha\in\A^s(L)\, :\, p^{-1}(\alpha) {\rm\ contains\ a\ good\ specialization} \}
$$ 
is not $L$-thin.

Since every point in $T$ has a pre-image in $W(\bar{L})$ of degree $\leq \deg(p)$, by Hilbert Irreducibility Theorem in the form of [Se, \S 9.2, Proposition 1], there exists an irreducible component of $W^\prime$ of $W$ on which the restriction of $p\circ \pi$ has degree $\leq\deg(p)$: this means precisely that $\pi$, when restricted to $W^\prime$,  has degree one, so it admits a rational section.

\cvd 
\medskip

Finally we remove the hypothesis that $U\subset\kappa^*$, proving Theorem 3.5 with $h=1$:
\medskip

\noindent{\bf Proposition 3.13}. {\it Suppose $\kappa$ is a finitely generated field, $U\subset\bar{\kappa}^*$ a torsion-free finitely generated multiplicative group, invariant under Galois conjugation. Let $f(T,n)\in\r_{\kappa,U}[T]$ be a monic polynomial. Suppose that for each $n\in\N$ the equation $f(t,n)=0$ has a solution in $\kappa$. Then there exists  an exponential polynomial $\t\in\r_{\kappa,U}$ such that identically $f(\t(n),d!n)\equiv 0$.
}
\smallskip

\noindent{\it Proof}. The deduction of Proposition 3.13 from 3.12 parallels the proof of Proposition 3.10. We only sketch the argument. We deduce from 3.12, applied with $\kappa(U)$ instead of $\kappa$, the existence of a functional solution  $\t\in\r_{\kappa(U),U}$; if such a solution is not in $\r_{\kappa,U}$, then it takes values in  $\kappa$ only for finitely many values of $n\in\N$. Then dividing the polynomial  $f(T,d!\cdot n)$ by the product of binomials $(T-\t^\sigma)$, where $\t^\sigma$ are the Galois conjugates of $\t$, we apply again Proposition 3.12 to the quotient polynomial. 
\cvd
\medskip

Before  proving the full Theorem 3.5, in several variables, it will be useful to give  a geometric formulation of the just proved Proposition 3.13:\medskip

\noindent{\bf Proposition  3.13 bis}. {\it Let $H$ be a commutative algebraic group, defined over $\kappa$, isomorphic to either a $\kappa$-torus $\T$ or to a product $\Ga\times\T$. Let $\Gamma\subset H(\kappa)$ be a Zariski-dense cyclic semigroup. Let $V$ be an affine algebraic variety, defined over $\kappa$, with each irreducible component of the same dimension as $H$; let $\pi:V\rightarrow H$ a finite map defined over $\kappa$ of degree $d\geq 1$. Suppose that  $\Gamma\subset\pi(V(\kappa))$. Then there exists an isogeny $\mu:H\rightarrow H$ sending  $g\mapsto g^{d!}$ (using multiplicative notation)  and a morphism $\theta:H\rightarrow V$, also defined over $\kappa$, such that $\pi\circ\theta=\mu$. 
}
\smallskip

\noindent{\it Proof}. First of all we can view $H$ as an algebraic subgroup of $\GL_N$, for some integer $N\geq 2$, defined over $\kappa$. So, its elements will be considered as $N\times N$ matrices. In particular the elements of $\Gamma$ will be matrices in $\GL_N(\kappa)$.

Let $V_1,\ldots,V_m$ be the irreducible components of $V$; for $j=1,\ldots,m$ let $\varphi_j\in\kappa[V]$ be a generator of the function field $\kappa(V_j)$ over $\pi^*(\kappa(H))$; since $\pi$ is by assumption a finite map,  $\varphi_j$ is integral over $\pi^*(\kappa[H])$. Put $V_j^\prime=\spec(\pi^*(\kappa[H])[\varphi_j])$ and $V^\prime=\spec(\kappa[V_1^\prime]\times\ldots\times\kappa[V_n^\prime])$ which corresponds to the disjoint union $V_1^\prime\cup\ldots\cup V_m^\prime$; it is an affine variety, birationally isomorphic to $V$. The minimal equation for $\varphi_j$ is given by a monic polynomial $f_j(T,g)\in\kappa[H][T]$; let $f(T,n)=f_1(T,n)\cdots f_m(T,n)\in\kappa[H][T]$ be the product of the minimal polynomials for $\varphi_1,\ldots,\varphi_m$. Let  $\gamma\in\Gamma$ be a generator of $\Gamma$; replacing if necessary $\gamma$ by some power of it, we can suppose that its spectrum generates a torsion-free subgroup in $\bar{\kappa}^*$ (and since $H$ is connected, every power of $\gamma$ still generates a Zariski dense semigroup). The polynomial $f(T,\gamma^n)$, which will be denoted for simplicity $f(T,n)$, will have its coefficients in a ring $\r_{\kappa,U}$, for a suitable torsion-free finitely generated group $U\subset\bar{\kappa}^*$, invariant under Galois conjugation. Proposition 3.13 assures the existence of a functional solution $\t\in\r_{\kappa,U}$ to the equation $f(\t,d! n)=0$. Consider the isogeny $\mu:H\rightarrow H$ sending $g\mapsto g^{d!}$. Via the identification $\kappa[H]\simeq\r_{\kappa,U}$, it acts on $\r_{\kappa,U}$ by sending the exponential polynomial $\y\in\r_{\kappa,U}$ to the exponential polynomial $\y^\prime$, where $\y^\prime(n):=\y(d!n)$ (see Remark 3.3). The functional solution $\t$ corresponds to a section $\theta:H\rightarrow V^\prime$ satisfying  $\pi\circ\theta=\mu$. Also, since $\kappa[H]$ is integrally closed, $\theta^*$ can be continued to $\kappa[V]$, so the morphism $\theta$   lifts to a regular map $H\rightarrow V$, obtaining the required section.
\cvd
\medskip

\noindent{\bf Lemma 3.14}. {\it Let $\kappa,U$ be as before; let $H=\Ga\times{\Bbb T}$ be a commutative algebraic group defined over $\kappa$ such that the ring $\r_{\kappa,U}$ is isomorphic to the $\kappa$-algebra $\kappa[H]$. Let $Z$ be an affine algebraic variety, also defined over $\kappa$, and let $f(T,z,n)\in\kappa[Z]\otimes_\kappa\r_{\kappa,U}[T]$ be a polynomial, monic in $T$, with coefficients in the ring $\kappa[Z]\otimes_\kappa\r_{\kappa,U}\simeq\kappa[Z\times H]$. Let $\{\alpha_1,\ldots,\alpha_l\}$ be a Galois invariant subset of $U$. 

There exists an affine algebraic variety $W$, defined over $\kappa$, and a morphism $p:W\rightarrow Z$, of finite degree, with the following property: for each $z_0\in Z(\bar{\kappa})$, there exists a functional solution $\t\in\r_{\bar{\kappa},U}$, whose roots belong to the set $\{\alpha_1,\ldots,\alpha_l\}$, to the equation
$$
f(t,z_0,n)\equiv 0\eqno(3.6)
$$
if and only if the fiber $p^{-1}(z_0)$ is not empty. Also, if $z_0$ is $\kappa$-rational, the (possible) functional solutions $\t\in\r_{\kappa,U}$ to the above equation are in bijection with the rational points of the fiber $p^{-1}(z_0)$.
}
\smallskip

In another language, the polynomial $f(T,z,n)$ defines a hypersurface $V\subset{\Bbb A}^1\times Z\times H$, which is naturally endowed with  projections $\pi_1: V\rightarrow Z$ and $\pi_2:V\rightarrow H$, so also  $\pi=(\pi_1,\pi_2) :V\rightarrow Z\times H$. The Lemma says that the (possible) sections of $\pi$ are parametrized by a $\kappa$-variety $W$, endowed with a projection $p:W\rightarrow Z$, whose rational points correspond to sections defined over $\kappa$. More precisely, the possible sections of $\pi_2$ over $\pi_1^{-1}(z)$ correspond to points of $p^{-1}(z)$. 
\smallskip

\noindent{\it Proof}. First of all, we know that the degree in $n$ of any possible functional solution is bounded by the degree in $n$ of the polynomial coefficients of $f(T,z,n)$; let $K$ be such a bound. 

Consider first the case where $U\subset\kappa^*$. In this case we argue as in the proof of Proposition 3.12;  write an unknown functional solution as
$$
\t(n)=\t(z,n)=p_1(z,n)\alpha_1^n+\ldots+p_l(z,n)\alpha_l^n\eqno(3.7)
$$
where each polynomial $p_j(X,z)\in\kappa[Z][X]$ will be of the form 
$$
p_j(X,z)=a_{0,j}(z)X^K+\ldots+a_{K,j}(z).\eqno(3.8) 
$$
The condition that for a given $z_0\in Z(\bar{\kappa})$ the exponential polynomial (3.7) be a functional solution to (3.6) amounts to an algebraic condition on the unknown coefficients $a_{0,j}(z_0),\ldots,a_{K,j}(z_0)$ of the polynomial $p_j$. Also, since we are supposing for the moment that all the roots $\alpha_1,\ldots,\alpha_l$ are $\kappa$-rational, the exponential polynomial (3.7) will belong to the ring $\r_{\kappa,U}$ (i.e. will be fixed by Galois conjugation over $\kappa$) if and only if the $l(K+1)$ coefficients $a_{i,j}$ are all $\kappa$-rational.
Hence the requested variety $W$ will be simply constructed as the closed subset of ${\Bbb A}^{l(K+1)}\times Z$ formed of the pairs 
$$
\{((a_{i,j})_{i,j},z)\in  {\Bbb A}^{l(K+1)}\times Z\, :\, \ 
{\rm s. t.\ the\ exponential\ polynomial\ } (3.7)\ {\rm satisfies}\ (3.6)\}
$$
The projection $p:W\rightarrow Z$ will be the natural projection on the second factor; it is a morphism of finite degree since for each choice of $z_0\in Z(\bar{\kappa})$, the specialized polynomial $f(T,z_0,n)$ will admit at most $\deg_T(f)$ functional solutions in $\r_{\bar{\kappa},U}$.

In the general case, when not all roots are rational, we must argue slightly differently. 
For each index $0\leq i\leq K$, write the (unknown) exponential polynomial $\t(z,n)$ as
$$
\t(z,n)=\t_0(z,n)+n\t_1(z,n)+n^2\t_2(z,n)+\ldots+n^K\t_K(z,n)
$$
for suitable (unknown) exponential polynomials $\t_0,\ldots,\t_K$ of the form
$$
\t_i(z,n)=a_{i,1}(z)\alpha_1^n+\ldots+a_{i,l}(z)\alpha_l^n.\eqno(3.9)
$$
Then $\t(z,n)$ is fixed by Galois conjugation if and only if each $\t_i$ has this same property. 
Now notice that each exponential polynomial $\t_i(z,n)$ always satisfies a linear recurrence relation defined over $\kappa$, of order $t$  (see [vdP, \S 2]); by this we mean that such a recurrence holds independently of the unknown coefficients $a_{i,j}$. Hence, by Lemma 3.1, in order that $\t_i(n,z)$ be fixed for the Galois action over $\kappa$, it is necessary and sufficient that it takes $\kappa$-rational values at the $l$ consecutive points $n=0,1,\ldots, l-1$. Then define, for each $i=0,\ldots,K$, the vector $(b_{i,1},\ldots,b_{i,t})=(b_{i,0}(z),\ldots,b_{i,t-1}(z))$ by putting
$$
b_{i,j}=\sum_{m=1}^l  \alpha_{m}^{j-1}a_{i,m}.
$$
Since the Van der Monde matrix $(\alpha_{m}^{j-1})_{1\leq l,j\leq l}$ is non-singular, the coefficients $a_{i,j}$ can be recovered from the $b_{i,j}$; hence the exponential polynomial (3.7) is a function of the $b_{i,j}$. Now, an exponential polynomial of the form (3.9) will be fixed by Galois conjugation if and only if $b_{i,j}$  all lie in $\kappa$. We then define $W\subset\A^{t(K+1)}\times Z$ as the closed set formed by the points $(b_{i,j},z)$ such that the exponential polynomial (3.7) satifies (3.6); as before this fact corresponds to an algebraic equation defining an affine $\kappa$-algebraic variety $W$.
\cvd
\medskip

\noindent{\it Proof of Theorem 3.5}. We argue by induction on $h$; the case $h=1$ is just Proposition 3.13. Then suppose Theorem 3.5 holds for exponential polynomials in $h-1$ variables. Our argument mimics the one in the proof of Proposition 3.13. Namely, we bound {\it a priori} the possible roots (and the degree) of a functional solution $\t$ with respect to the last variable $n_h$. Again, we fix a basis for the torsion-free multiplicative group $U$, so that we can speak of the height of any element in $U$. The crucial observation is that, for every given vector $(n_1,\ldots,n_{h-1})\in\N^{h-1}$, the roots of the exponential polynomial $(n\mapsto \y_i(n_1,\ldots,n_{h-1},n))$ belong to a fixed set, independent of $(n_1,\ldots,n_{h-1})$. So, by Lemma 3.8, the height of any possible functional solution to (3.5) is bounded independently of $(n_1,\ldots,n_{h-1})$; the same holds for the degree in $n_h$ of its polynomial coefficients. Let $\{\alpha_1,\ldots,\alpha_l\}$ be the corresponding set of possible roots; up to adding, if necessary, their conjugates, we can suppose that the finite set $\{\alpha_1,\ldots,\alpha_l\}$ is invariant by Galois conjugation over $\kappa$. Let $K$ be a bound for the degree of the polynomial coefficients of $\t$ in $n_h$ (for any possible functional solution). We search for exponential polynomials of the form
$$
\t(n)=\t_{(n_1,\ldots,n_{h-1})}(n)=p_1(n)\alpha_1^n+\ldots+p_l(n)\alpha_l^n
$$
satisfying $f(\t(n),d!n)\equiv 0$.
Here the polynomials $p_1(X),\ldots,p_l(X)\in\bar{\kappa}[X]$ have degree $\leq K$, and their coefficients depend on $(n_1,\ldots,n_{h-1})$.
We apply Lemma 3.13, with $Z=H^{h-1}$, constructing an affine $\kappa$-algebraic variety $W$ and a morphism $p:W\rightarrow H^{h-1}$ with the property of Lemma 3.14. 

 Then, up to changing $W$ by a variety birationally equivalent to it, we can suppose $W\subset\A^1\times H^{h-1}$ is given by a single monic equation of the form $g(T,n_1,\ldots,n_{h-1})=0$, where $g(T,n_1,\ldots,n_{h-1})\in\r_{\kappa,U}^{\otimes (h-1)}[T]$. The  hypothesis of Theorem 3.5 implies that for every $(n_1,\ldots,n_{h-1})\in\N^{h-1}$, there exists a rational solution $s\in\kappa$  to the equation 
$$
g(s,n_1,\ldots,n_{h-1})=0.
$$
The inductive hypothesis implies the existence of a  functional solution ${\goth s}\in\r_{\kappa,U}^{\otimes (h-1)}$, for $(n_1,\ldots,n_{h-1})\in d!\cdot \Z^{h-1}$. This corresponds to a section for $p$: in other words, the unknown coefficients of $p_1,\ldots,p_l$ can be written as exponential polynomials in $h-1$ variables $(n-1,\ldots,n_{h-1})$, which implies that $\t$ can be written as an exponential polynomial in $h$ variables, as wanted.
\cvd

\bigskip

\noindent {\bf \S 4. Auxiliary results on linear algebraic groups.} 

We collect and prove in this section some lemmas of geometric nature to be used in the proof of our main theorems. As in the previous sections, $\kappa$ will be a fixed finitely generated field imbedded in the field $\C$ of complex numbers. \medskip

\noindent{\bf Lemma 4.1}. {\it Let $\Gamma\subset\GL_N(\C)$ be a semigroup of invertible matrices. Then its Zariski closure in $\GL_N$ is an algebraic subgroup of $\GL_N$. } \medskip

This is already known; a proof can be found for instance in [Be, \S 1]. We give here an alternative proof, closer in spirit to the techniques involved in the present paper. \smallskip

\noindent{\it Proof}. Let $\bar{\Gamma}$ be the Zariski closure of $\Gamma$ in $\GL_N$. Clearly it is a semigroup; to prove it is a group, it suffices to show that it is closed under the map $g\mapsto g^{-1}$. Let $g\in\bar{\Gamma}$; then its positive powers all belong to $\bar{\Gamma}$; we now show that under this hypothesis its negative powers too belong to $\bar{\Gamma}$. This amounts to showing that every regular function on $\GL_N$, vanishing on the positive powers of $g$, also vanishes on the negative ones. Now such a function can be written as a polynomials function of the entries of the matrices of $\GL_N$, possibly divided by the determinant function. Since the entries of the powers $g^n$ of $g$ are linear recurrent sequences in $n$, every regular function on $\GL_N$, calculated in the sequence $g^n$, is a linear recurrent sequence of the variable $n$. Hence, if it vanishes for all positive $n$, it also vanishes for negative $n$.
\cvd

\bigskip

\noindent{\bf Lemma 4.2}. {\it Let $G$ be a linear algebraic group, $H\subset G$ an algebraic subgroup, both defined over the finitely generated field $\kappa$.  There exists an integer $m=m(G,H,\kappa)>0$ with the following property: for every matrix $g\in G(\kappa)$ such that some positive power of $g$ belongs to $H$, we have $g^m\in H$. } \medskip

\noindent{\it Proof}. By [Bo 2, Theorem 5.1] $H$ is the stabilizer of a point $P\in\Pr_M(\kappa)$ for a suitable immersive linear representation of $G$ in $\GL_{M+1}$ defined over $\kappa$. Then, considering $G$ imbedded in $\GL_{M+1}$ and acting canonically on the projective space $\Pr_M$, each element $g\in G$ such that $g^n\in H$ for an integer $n$ corresponds to a matrix $g$ such that $g^n$ fixes the point $P$. Let $n=n(g)$ be the minimal positive integer with such property and suppose it is $>1$. Then $g$ does not fix $P$, nor does any power $g^m$ with $1<m<n$,  but its $n$-th power $g^n$ does fix $P$. So $g$ admits two eigenvalues whose ratio is a primitive $n$-th root of $1$. If $g$ is defined over $\kappa$, this implies that the $n$-th roots of unity have degree $\leq M+1$ over $\kappa$, and this facts gives a bound on $n$, since $\kappa$ is finitely generated. Letting $m$ be the least common multiple of the possible values of $n$ we obtain an integer with the sought property.\cvd

\medskip

Given an algebraic subvariety $Z\subset G$ of an algebraic group $G$ and an element $g\in G$, we shall denote by $Z\cdot g$ the image of $Z$ under the right-translation by $g$. The following result is a generalization of Schur's Theorem on finitely generated matrix groups [CR, Theorem 36.2]; indeed Schur's Theorem can be easily recovered from the case $Z=\{1_G\}$ below: \medskip

\noindent{\bf Lemma 4.3}. {\it Let $G$ be a linear algebraic group defined over the finitely generated field $\kappa$, $Z\subset G$ be a closed algebraic subvariety defined over $\kappa$ of strictly inferior dimension. Let $\Gamma\subset G(\kappa)$ be a semigroup, Zariski-dense in $G$. There exists an element $\gamma\in \Gamma$ such that the algebraic group generated by $\gamma$ is connected and no positive power $\gamma^n$ of $\gamma$ satisfies $Z\cdot\gamma^n=Z$. } \medskip

\noindent{\it Proof.} Let $H$ be the subgroup of $G$ formed by the elements $g\in G$ such that $Z\cdot g=Z$. Clearly $H$ is an algebraic subgroup of $G$. Moreover $H$ has dimension $<\dim(G)$ if $Z$ satisfies the same condition (which we assumed). Let us first prove the existence of an element $\gamma\in\Gamma$ such that $Z\cdot\gamma^n\neq Z$ for all $n>0$. Thanks to the previous lemma, it suffices to prove that for every $m\geq 1$ there exists an element $\gamma\in\Gamma$ such that $\gamma^m\not\in H$. Assume the contrary, so $\gamma^m$ does belong to $H$ for some fixed $m\geq 1$ and all $\gamma$. Hence $\Gamma$ would be contained in the algebraic subvariety of $G$ defined by the condition $g^m\in H$; since $\Gamma$ is Zariski dense, $G$ itself would be contained in such a variety, so every element of $g$ would satisfy such a relation. To see that this is impossible, just consider a one-parameter subgroup of $G$ not contained in $H$, which exists due to the hypothesis $\dim(H)<\dim(G)$. Such a subgroup contains only finitely many elements whose $m$-th powers lie  in $H$, obtaining a contradiction. Replacing $\gamma$ by a suitable power, the previous condition still holds and we obtain moreover an element generating a connected algebraic group.\cvd

\medskip

All the previous lemmas were meant to prove the following Proposition, for which we introduce one more definition: for a $h$-tuple $(g_1,\ldots,g_h)\in G(\C)^h$, $G$ being as usual an algebraic group, we let $\s(g_1,\ldots,g_h)\subset G(\C)$ be the set 
$$ 
{\cal S}(g_1,\ldots,g_h):= \{g_1^{n_1}\cdots g_h^{n_h}\, :\, (n_1,\ldots,n_h)\in\N^h\}.\eqno(4.1) 
$$
 \medskip

\noindent{\bf Proposition 4.4}. {\it Let $G$ be a connected linear algebraic group defined over the finitely generated field $\kappa$, $\Gamma\subset G(\kappa)$ a Zariski-dense semigroup. There exist an integer $h\geq 1$ and elements $\gamma_1,\ldots,\gamma_h\in\Gamma$ with the following property: for every positive integer $D$ the set $\s(\gamma_1^D,\ldots,\gamma_h^D)$ is Zariski-dense in $G$. Moreover, letting, for each index $i=1,\ldots,h$, $Z_i$ denote the Zariski closure of $\s(\gamma_1^D,\ldots,\gamma_i^D)$, the algebraic varieties $Z_i$ are irreducile and pairwise distinct (and $Z_h=G$). Also, for each $i=1,\ldots,h$, the algebraic subgroup generated by $\gamma_i^D$ is connected. } \medskip

Note that Proposition 4.4 fails if one omits the arithmetic condition on the field $\kappa$. For instance if $G=\Gm$ and $\Gamma\subset \Gm(\C)$ is the torsion subgroup of $G$, i.e. the group of roots of unity, then $\Gamma$ is Zariski dense in $G$ but every set of the form (4.1), with $g_i\in\Gamma$,  is finite. \medskip

\noindent{\it Proof of Proposition 4.4}. If $G$ is zero dimensional, then $\Gamma=G=\{1_G\}$ and we are done. Suppose $G$ has positive dimension. By the Lemma 4.3, applied with $Z=\{1_G\}$, there exists an element $\gamma_1$ of infinite order. Then the Zariski closure $Z_1$ of the cyclic group generated by $\gamma_1$ has positive dimension. Replacing $\gamma_1$ by a suitable power of it if necessary, we can suppose that $Z_1$ is irreducible. If the dimension of $Z_1$ coincides with the dimension of $G$ we are done, since $G$ is connected. Otherwise, by the preceding lemma, there exists an element $\gamma_2\in\Gamma$ such that the sets $Z\cdot\gamma_2^n$ for $n=1,2,\ldots$ are pairwise distinct. Again by taking a suitable power of $\gamma_2$, we can suppose it generates a connected algebraic group $H_2$. Then the Zariski closure $Z_2$ of the set $\s_2=\s(\gamma_1,\gamma_2)$ has strictly larger dimension (since it contains infinitely many pairwise distinct subvarieties isomorphic to $Z_1$) and is connected (since it is the image of the connected variety $Z_1\times H_2$ under the map $(g,h)\mapsto g\cdot h$). Note that $Z_2$, unlike $Z_1$, need  not be an algebraic subgroup. After at most $\dim(G)$ steps we reach a variety $Z_h$ which has the same dimension as $G$, hence coincides with $G$; we also obtain a sequence $H_1,\ldots,H_h$ of connected algebraic groups of positive dimension generated by $\gamma_1,\ldots,\gamma_h$ respectively. We now observe that given such $\gamma_1,\ldots,\gamma_h$, the varieties $Z_1,\ldots, Z_h$ and the subgroups $H_1=Z_1,\ldots, H_h$ do not change if we replace the $\gamma_i$ by any positive power $\gamma_i^D$; this follows from the fact that if the connected algebraic group $H_i$ is generated by an element $\gamma_i$, it is also generated by any positive power $\gamma_i^D$. \cvd

\medskip

We shall later choose the integer $D$ in the above Lemma in such a way that the spectra of the matrices $\gamma_1^D,\ldots,\gamma_h^D$ generate  a torsion free group (compare with the notion of ``sous-groupe net" in [Bo 1]). We shall prove \medskip

\noindent{\bf Proposition 4.5}. {\it Let $G$ be a connected linear algebraic group defined over the finitely generated field $\kappa$, $\Gamma\subset G(\kappa)$ a Zariski-dense sub-semigroup of $G$. There exist an integer $h\geq 1$ and elements $\gamma_1,\ldots,\gamma_h\in\Gamma$ with the following properties: letting, for $i=1,\ldots,h$, $H_i$ be the Zariski closure of the group generated by $\gamma_i$, 

\item{(i)} the algebraic groups $H_i$ are irreducible; 

\item{(ii)} the map $H_1\times\ldots\times H_h\rightarrow G$ sending $(x_1,\ldots,x_h)\mapsto x_1\cdots x_h$ is surjective;

\item{(iii)} the spectra of the matrices $\gamma_1,\ldots,\gamma_h$, generate in $\C^*$ a torsion free subgroup. }
\smallskip

\noindent{\it Proof}. 
Having choosen $\gamma_1^\prime,\ldots,\gamma_h^\prime$ as in Proposition 4.4 we let $D$ be the order of the torsion subgroup of the multiplicative group generated by the spectra of $\gamma_1^\prime,\ldots,\gamma_r^\prime$. Then putting $\gamma_i:=\gamma_i^{\prime D}$ we obtain the all the conditions {\it (i), (ii), (iii)} are satisfied.
\cvd

\medskip

We are interested in the  Stein factorization of the surjective map $H_1\times\ldots\times H_h\rightarrow G$ appearing in Proposition 4.5. In the sequel of this paragraph, $\kappa$ denotes any field.
\smallskip

Put $H:=H_1\times\ldots\times H_h$ and let $\psi:H\rightarrow G$ be the surjective map  
$$
\psi:\, H\ni(x_1,\ldots,x_h)\mapsto \psi(x_1,\ldots,x_h)=x_1\cdots x_h\in G.
$$ 

\noindent{\bf Definition}. We say that an automorphism $\sigma$ of $H$ (in the sense of $\kappa$-algebraic varieties) {\it preserves the fibers} of $\psi$ if there exists an automorphism $\bar{\sigma}$ of $G$ (as algebraic variety) such that $\psi\circ\sigma=\bar{\sigma}\circ\psi$. Clearly, such an automorphism is uniquely determined by $\sigma$.
\smallskip

The fiber preserving automorphisms of $H$ form a group; they are characterised by the following property: for every choice of points $\alpha_1,\alpha_2\in H$, $\psi(\alpha_1)=\psi(\alpha_2)$ if and only if $\psi(\sigma(\alpha_1))=\psi(\sigma(\alpha_2))$.
\medskip

The next lemma garantees that the group of automorphisms of $G$ of the form $\bar{\sigma}$, for some $\sigma:H\rightarrow H$ preserving the fibers of $\psi$, acts transitively on $G$. 
\medskip

\noindent{\bf Lemma 4.6}. {\it Let $G$ be a connected linear algebraic group defined over a field $\kappa$, $H_1,\ldots,H_h$ be connected subgroups, also defined over $\kappa$. Let $\psi:H_1\times \ldots\times H_h=H\rightarrow G$ be the map defined above, and suppose it is surjective.
Then the fibers of $\psi$ are all isomorphic. Also, given two points $g_1,g_2$ of $G$, there exists an automorphism of $H$,  preserving the fibers of $\psi$, and sending $\psi^{-1}(g_1)$ to $\psi^{-1}(g_2)$.
}
\smallskip

\noindent{\it Proof}. Let $g\in G$ be a point. Write $g=a_1\cdot a_2\cdots a_h$, so that $\psi^{-1}(g)$ contains the point $(a_1,a_2,\ldots,a_h)$. We first prove that the fiber of $g$ is isomorphic to the fiber of $g^\prime:=a_2\cdots a_h$, and that such an isomorphism can be chose to preserve all fibers: clearly the automorphism $\sigma$ of $H$: $\sigma(x_1,\ldots,x_h)\mapsto (a_1^{-1}x_1,\ldots,x_h)$ preserves the fibers and defines by restriction to $\psi^{-1}(g)$  an isomorphism between $\psi^{-1}(g)$ and $\psi^{-1}(g^\prime)$. (Note that $\bar{\sigma}:G\rightarrow G$ is the left-translation by $a_1^{-1}$).
For the same reason, there is a fiber preserving automorphism  sending the fiber of $g^\prime$ to the fiber of $g^\second:=a_3\cdots a_h$. After $h-1$ steps we obtain that the fiber of $g$ is isomorphic, via a fiber preserving automorphism of $H$, to the fiber of the neutral element $1_G$. By transitivity, we obtain the Lemma.
\cvd
\medskip

\noindent{\bf Proposition 4.7}. {\it Let the subgroups $H_1,\ldots,H_h$ and the map $\psi$ be as in Lemma 4.6. Put $H:=H_1\times\ldots\times H_h$. Then there exists a connected algebraic group $G^\prime$, defined over $\kappa$, such that the map $\psi:H\rightarrow G$ factors as $\psi=\psi_1\circ \psi_2$, where $\psi_2:H\rightarrow G^\prime$ has irreducible fibers and $\psi_1:G^\prime\rightarrow G$ is an isogeny.
}
\medskip

The above factorization of $\psi$ will be referred to as the {\it Stein factorization}. 
\medskip

\noindent{\it Proof of Proposition 4.7}. The surjective morphism $\psi:H\rightarrow G$ induces an injection $\psi^*: \kappa[G]\hookrightarrow\kappa[H]$ of $\kappa$-algebras; let $G^\prime$ be the affine variety corresponding to the integral closure of $\psi^*(\kappa[G])$ in the function field $\kappa(H)$. Since $H$ is a normal variety, it is a subring of $\kappa[H]$. We then have a factorization $H\rightarrow G^\prime\rightarrow G$ of the map $\psi$, where the second arrow, $\psi_1:G^\prime\rightarrow G$, is a finite map.

  We first prove that the variety $G^\prime$ is smooth. Let us first observe that for the fiber preserving automorphisms $\sigma$ of $H$ induce automorphisms of $G^\prime$; in fact, from $\psi\circ\sigma=\bar{\sigma}\circ\psi$, it follows that $\sigma^*$ preserves the subring $\psi^*(\kappa[G])$, so also its integral closure $\psi_1^*(\kappa[G^\prime])$. We shall denote by  $\sigma^\prime$ the induced automorphism of $G^\prime$, which satisfies $\sigma^\prime\circ\psi_1=\psi_1\circ\sigma$. 
Let $Z$ be the singular locus of $G^\prime$ and suppose by contradiction that it is non empty. We shall prove that $\psi_1(Z)=G$, which is impossible since $\dim(Z)<\dim (G^\prime)=\dim(G)$.
Suppose then that $z\in Z$ is a singular point of $G$; let $g\in G$ be any point. Choose a pre-image $a\in \psi_2^{-1}(z)$ for $z$ and let $\sigma$ be a fiber preserving automorphism of $H$ sending $z$ to a point in $\psi^{-1}(g)$.  Then the induced automorphism $\sigma^\prime$ of $G^\prime$ sends $z$ to a point $\sigma^\prime(z)$ such that $\psi_1(\sigma^\prime(z))=g$. Since $\sigma^\prime$ is an automorphism, $\sigma^\prime(z)$ is also a singular point, so it belongs to $Z$, thus proving that $g\in\psi_1(Z)$.

  We now prove that the finite map $\psi_1:G^\prime\rightarrow G$ is an unramified cover, i.e. its differential is surjective at every point of $G^\prime$.  
The proof is very similiar to the previous one. Let $W\subset G^\prime$ be the (closed, proper) subset of points where the differential of $\psi_1$ is not surjective and suppose by contradiction it is non empty. We shall prove that $\psi_1(W)=G$, obtaining a contradiction, since at each point of the variety $W$ the differential of $\psi_1$ has rank $<\dim(G)$. 
Let as before $g\in G$ be any point of $G$ and $w\in W$ be a point of $W$. As before, there exists an automorphism  $\sigma^\prime$ of $G^\prime$, induced by a fiber preserving automorphism $\sigma$ of $H$, satisfying $\psi_1(\sigma^\prime(w))=g$. This proves that $g\in\psi_1(W)$, obtaining the sought contradiction.

 It remains to prove the second part of the Proposition, namely the irreducibility of the fibers of $\psi_2$. The generic fiber of  $\psi_2:H\rightarrow G^\prime$ is irreducible, since $\psi_2^*(\kappa(G^\prime))$ is algebraically closed in $\kappa(H)$.
Every point of $G$ has $\deg(\psi_1)$ preimages in $G^\prime(\bar{\kappa})$. Hence there exists a point $g\in G$ such that $\psi^{-1}(g)$ has $\deg(\psi_1)$ preimages in $G^\prime$, and each of such preimages has irreducible fiber with respect to $\psi_2$ (actually an open dense subset of $G$ of such points with this property). 
Then the number of irreducible components of $\psi^{-1}(g)$ of such a point $g$ equals the degree of $\psi_1$. We shall prove that it is so for every other point of $G$.
Now, since the variety $\psi^{-1}(g)=\psi_2^{-1}(\psi_1^{-1}(g))$ is isomorphic to $\psi^{-1}(x)=\psi_2^{-1}(\psi_1^{-1}(x))$, the fibers have the same number of irreducible components. But now, since    $\psi_1^{-1}(x), \psi_1^{-1}(g)$ have the same cardinality  $\deg(\psi_1)$,  each fiber, with respect to $\psi_2$, of each point in $\psi_1^{-1}(x)$ must be irreducible, concluding the proof.
\cvd

We end this section with a result of different nature, which will be used in the proof of Theorem 1.2: it is probably well known, and can be proved in several different ways, but we cannot locate any reference to it on the literature. The proof given below follows a suggestion of U. Zannier.
\medskip

\noindent{\bf Proposition 4.9}. {\it Let $G$ be a connected linear algebraic group, $\chi\in\kappa[G]$ be a regular function. Suppose that $\chi$ has no zeros in $G(\bar{\kappa})$ and that $\chi(1_G)=1$. Then $\chi$ is a character.
}
\smallskip

\noindent{\it Proof}. Let $\overline{G}$ be a normal compactification of $G$, i.e. a normal complete variety containing $G$ as an open subset (in the Zariski topology). Let $D_1,\ldots,D_t$ be the component of the divisor at infinity, i.e. of the hypersurface $\bar{G}\setminus G$. For each $g\in G(\bar{\kappa})$, denote by $\chi^g$ the regular function $x\mapsto \chi^g(x)=\chi(gx)$. All the zeros and poles of $\chi^g$, if any, are contained in the hypersurfaces $D_1,\ldots,D_t$. Hence, for each $g\in G$, the divisor of poles of the rational function $\chi^g\in\kappa(G)$ can be written as
$$
(\chi^g)_\infty=a_1(g)D_1+\ldots+a_t(g)D_t=:A(g)
$$
for suitable non-negative integers $a_1(g),\ldots,a_t(g)$. The same holds for the divisor of zeros, namely
$$
(\chi^g)_0=b_1(g)D_1+\ldots+b_t(g)D_g=:B(g),
$$
for non-negative integers $b_1(g),\ldots,b_t(g)$. Since $G$ is an irreducible variety, the effective divisors $A(g)$, for $g\in G$, are all algebraically equivalent, so in particular they have the same degree with respect to a projective imbedding $\bar{G}\hookrightarrow \Pr_N$ ; the same is true of the divisors $B(g)$. So the coefficients $a_1(g),\ldots,a_t(g),b_1(g),\ldots,b_t(g)$ have only finitely many possibilities. Hence there exists a Zariski open set $\Omega\subset G$ such that for every $g\in\Omega$ the (principal) divisor $(\chi^g)=B(g)-A(g)$ is fixed,   equal to the divisor of a  function $\chi^{g_0}$, for a fixed $g_0\in G$. Then for all $g\in \Omega$ the ratio $\chi^g/\chi^{g_0}$ is a constant function on $G$, say $\rho(g)$. But the equality $\chi(gh)=\rho(g)\chi(g_0h)$, valid for all $(g,h)\in\Omega\times G$, immediately implies the same equality for all $(g,h)\in G\times G$. This, and the fact that $\chi(1_G)=1$, easily implies that $\chi$ is a character.
\cvd
\bigskip

\noindent {\bf \S 5. Proof of the main theorems}. 

Recall that $\kappa\subset\C$ always stands for a fixed finitely generated field of characteristic zero.

\medskip

We now begin the {Proof of Theorem 1.6}, which is the crucial ingredient in the proofs of Theorems 1.1 and 1.5. We shall first consider a slightly different statement, which will be proved to be equivalent to Theorem 1.6: 
\smallskip

\noindent {\bf Proposition 5.1}. {\it Let $V$ be an affine variety, with $\dim(V)=\dim(G)$.
Suppose that each irreducible component of $V$ has the same dimension, equal to the dimension of $G$. Let $\pi:V\rightarrow G$ be a  morphism, such that its restriction to each irreducible component of $V$ is dominant; let $\Gamma\subset G(\kappa)$ be a Zariski-dense semigroup with $\Gamma\subset\pi(V(\kappa))$. Then there exists an algebraic group $\tilde{G}$, an isogeny $p:\tilde{G}\rightarrow G$ and a rational map $\theta:\tilde{G}\rightarrow V$, all defined over $\kappa$, such that $\pi\circ\theta=p$.
}
\smallskip

\noindent {\it Proof}. Denote by $V_1,\ldots,V_m$ the irreducible components of $V$. The morphisms $\pi_{|V_j}:V_j\rightarrow G$ correspond to  inclusions of $\kappa$-algebras $\pi_{|V_j}^*:\kappa[G]\hookrightarrow\kappa[V_j]$. For each $j\in\{1,\ldots,m\}$, let $t_j\in\kappa[V_j]$ be a generator for the field extension $\kappa(V_j)/\pi_{|V_j}^*(\kappa(G))$, which is integral over $\pi_{|V_j}^*(\kappa[G])$. 
Putting $W_j={\rm spec}(\pi_{|V_j}^*(\kappa[G])[t_j])$ we obtain affine varieties $W_1,\ldots,W_m$, endowed with finite mappings  $\pi_j:W_j\rightarrow G$ and birational morphisms $\omega_j:V_j\rightarrow W_j$, with $\pi_j\circ\omega_j=\pi_{|V_j}$. The hypothesis that $\Gamma\subset\pi(V(\kappa))$ implies that 
$$
\Gamma\subset \bigcup_{j=1}^m\pi_j(W_j(\kappa)).
$$
Letting $W$ be the disjoint union of $W_1,\ldots,W_m$, so  $W:=W_1\cup\ldots\cup W_m={\rm spec}(\kappa[W_1]\times\ldots\times\kappa[W_m])$, we 
 are then reduced to a finite mapping, denoted again by  $\pi:W\rightarrow G$ with moreover the property that on each irreducible component $W_j$, the ring extension $\kappa[W_j]$ of $\pi_j^*(\kappa[G])$ is generated by a single element $t_j$. Of course, if we prove that on a suitable unramified (connected) cover $\tilde{G}\rightarrow G$ there exists a regular section for the map $W\rightarrow G$, we can deduce the same conclusion for $V$, which is birationally isomorphic to $W$, up to the fact that the section to $V$ would be only rational (not necessarely regular). Since this is exactly our thesis, we shall be content to prove the existence of a section $\tilde{G}\rightarrow W$. 

For each $j\in\{1,\ldots,m\}$, let  the minimal equation for $t_j$ be
$$
T^{d_j}+\varphi_{1,j}(g)T^{d_j-1}+\ldots+\varphi_{d_j,j}(g)=0.
$$
Here $d_j=\deg(\pi_{|V_j})=\deg(\pi_j)$ and $\varphi_{j,1},\ldots,\varphi_{j,d_j}$ are regular functions on $G$, identified via $\pi^*$ to regular functions on $V$. 

The hypothesis that $\Gamma \subset\pi(V(\kappa))$ implies that for each $g\in\Gamma$, there exists at least an index $j\in\{1,\ldots,m\}$ and a rational specialization of $t_j\mapsto t_j(\gamma)\in\kappa$ satisfying the above displayed  equation.
Taking the product of the polynomials 
$T^{d_j}+\varphi_{1,j}(g)T^{d_j-1}+\ldots+\varphi_{d_j,j}(g)\in\kappa[G][T]$, for $j=1,\ldots,n$, we obtain a polynomial $f(T,g)\in\kappa[G]$, of degree $d:=d_1+\ldots+d_m$, which we write as
$$
f(T,g)=T^d+\varphi_1(g)T^{d-1}+\ldots+\varphi_d(g)
$$
for suitable regular functions $\varphi_1,\ldots,\varphi_d\in\kappa[G]$.
By hypothesis, such a polynomial has a root in $\kappa$ for each specialization $g\in\Gamma$.

Let now $\gamma_1,\ldots,\gamma_h$, as well as $H_1,\ldots,H_h$, $H:=H_1\times\ldots\times H_h$, be as in Propositions 4.5, 4.6, 4.7. Denote by $\psi:H\rightarrow G$, as in Proposition 4.6, the map sending the $h$-tuple $H\ni(g_1,\ldots,g_h)\mapsto \psi(g_1,\ldots,g_h)=g_1\cdots g_h$.
Putting, for $i=1,\ldots,d$, 
$$
\y_i(n_1,\ldots,n_h)=\varphi_i(\gamma_1^{n_1}\cdots\gamma_h^{n_h})
$$
we obtain $d$ exponential polynomials in $\r_{\kappa,U}^{\otimes h}$, where $U$ is the torsion-free group generated by the spectra of $\gamma_1,\ldots,\gamma_h$. We remark at once that each $\kappa$-algebra $\kappa[H_i]$ imbeds canonically into $\r_{\kappa,U}^{\otimes h}$.  
Apply Theorem 3.5 to the polynomial $f(T,n_1,\ldots,n_h)\in\r_{\kappa,U}^{\otimes h}[T]$ defined by
$$
\eqalign{
f(T,n_1,\ldots,n_h)&:=T^d+\y_1(n_1,\ldots,n_h)T^{d-1}+\ldots+\y_d(n_1,\ldots,n_h)\cr &=
T^d+\varphi_1(\gamma_1^{n_1}\cdots\gamma_h^{n_h})T^{d-1}+\ldots+
\varphi_d(\gamma_1^{n_1}\cdots\gamma_h^{n_h}).
}
$$
We obtain a functional solution in the ring $\r_{\kappa,U }^{\otimes h}$, after replacing if necessary $\gamma_1,\ldots,\gamma_h$ by $\gamma_1^{d!},\ldots,\gamma_h^{d!}$. This corresponds geometrically to an unramified covering $H^\prime$ of $H:=H_1\times\ldots\times H_r$:
$$
\eta: H^\prime\rightarrow H
$$
and a regular map $\lambda:H^\prime\rightarrow W$ with $\pi\circ\lambda=\psi\circ\eta$.

Consider the morphism $\psi\circ\eta:H^\prime\rightarrow G$. We want to investigate its Stein factorization. Recall (Proposition 4.7) that $\psi$ factors as $\psi=\psi_1\circ\psi_2$, where $\psi_2: H\rightarrow G^\prime$ has irreducible fibers and $\psi_1: G^\prime\rightarrow G$ is an isogeny (connected unramified covering). Let now $\tilde{G}$ be the affine variety corresponding to the integral closure of $(\psi_2\circ\eta)^*(\kappa[G^\prime])$ in $\kappa[H^\prime]$; we obtain a factorization of $\psi_2\circ\eta$ as $\psi_2\circ\eta=\eta_1\circ\eta_2$ where $\eta_1: \tilde{G}\rightarrow G^\prime$ has finite degree and $\eta_2: H^\prime\rightarrow \tilde{G}$ has connected fibers. 
Now, since $\eta: H^\prime\rightarrow H$ is unramified, the morphism $\eta_1: \tilde{G}\rightarrow G^\prime$ is also unramified, so $\tilde{G}$ admits an algebraic  group structure such that $\eta_1$ becomes an isogeny.
Let $p:\tilde{G}\rightarrow G$ be the composite $p=\psi_1\circ\eta_1$; it is an isogeny with respect to the above mentioned algebraic group structure on $\tilde{G}$.

Our next goal is to prove that the map $\lambda:H^\prime\rightarrow W$ factors as $\lambda=\theta\circ\eta_2$, for a morphism $\theta:\tilde{G}\rightarrow W$, so that we can obtain the commutative diagram:
 $$
\matrix{
W &\longleftarrow & \tilde{G} & \longleftarrow & H^\prime \cr
\downarrow & {}  &\downarrow & {} &\downarrow \cr
G  & \longleftarrow & G^\prime  & \longleftarrow & H \cr
}
$$
This amounts to saying that $\lambda$ is constant on every fiber of $\eta_2$; now this is clear since (1) the fibers of $\eta_2$ are connected and (2) on each fiber of $\eta_2$ the value of $\lambda$ has only finitely many possibilities, since it must belongs  to the fiber with respect to $\pi$ of a single point in $G$.
\cvd
 
\medskip

The above Proposition can be generalized to reducible varieties $V$ of mixed dimension. The crucial point is the following 
\medskip

\noindent{\bf Lemma 5.2}. {\it Let $G$ be a normal irreducible affine variety, $Y\subset G$ a proper irreducible closed subvariety, both defined over $\kappa$. There exists an affine variety $W$, with $\dim(W)=\dim(G)$, and a finite map $\pi:W\rightarrow G$ of degree $>1$ such that $Y(\kappa)\subset\pi(W(\kappa))$, $\pi$ is ramified (over a hypersurface in $G$). 
} 
\smallskip

We remark at once that our Lemma 5.2 implies in particular that, in the notation of [Se, chap 9, p. 121], thin sets of type 1 are also of type 2.  
\smallskip

\noindent{\it Proof}. By imbedding $Y$ in a hypersurface defined over $\kappa$ we can reduce to the case  $Y$ is a hypersurface in $G$. Let $f\in\kappa[G]$ be a regular function having a zero of multiplicity one in $Y$ (such a function exists since $G$ is normal). Put $W:={\rm spec}(\kappa[G][\sqrt{f}])$. We obtain a degree two cover $\pi:W\rightarrow G$, defined over $\kappa$, ramified over the hypersurface $Y$; note that it is an isomorphism  on $\pi^{-1}(Y)$. In particular $Y(\kappa)\subset \pi(W(\kappa))$.
\cvd

\medskip

\noindent{\it Proof of Theorem 1.6}. (We shall prove at the same time that  in Proposition 5.1 one can omit the hypothesis that all components of $V$ have the same dimension.)
 
Suppose the hypotheses of Theorem 1.6 hold. By Lemma 5.2, we can construct a variety $W$ from $V$ as follows: consider the irreducible components of $V$ of dimension $<\dim(G)=\dim(V)$; for each such component $V^\prime$, remark that $\pi(V^\prime)$ is contained in a proper closed subvariety of $G$; let $W^\prime$ be an irrreducible affine variety with $\dim(W^\prime)=\dim(G)$, endowed  with a map $\pi_{W^\prime}$ such that (1) $\pi_{W^\prime}(W^\prime(\kappa))\supset \pi(V^\prime(\kappa))$ and (2) $\pi_{W^\prime}$ is ramified; the existence of such varieties and maps is assured by Lemma 5.2.
Let $W$ be a normalization of the disjoint union of the irreducible components of maximal dimension of $V$ and the varieties $W^\prime$, obtained as explained from the components of lower dimension of $V$.
The new affine variety $W$ is endowed with a morphism $\pi_W:W\rightarrow G$, satisfying the hypotheses of Proposition 5.1; namely $\pi(W(k))\supset \Gamma$. Also, $\pi_W$ coincides with $\pi$ on the union of the irreducible components of $V$ of maximal dimension. By  Proposition 5.1, there exists an isogeny $p:\tilde{G}\rightarrow G$ and a rational map $\theta:\tilde{G}\rightarrow W$ with $\pi\circ\theta=p$. Consider the image $V^\prime=\theta(\tilde{G})$ of the map $\theta$; it is an irreducible component of $W$.
To end the proof, we must show that it is a component of $V$. Now,   since $\pi\circ\theta$ is unramified, $V^\prime$ cannot be one of the components of $W$ constructed via Lemma 5.2 from those of lower dimension in $V$, otherwise $\pi$ (hence $\pi\circ\theta$) would be ramified, at least over the smooth locus of the ramification divisor of $\pi$. So $V^\prime$ is a component of $V$.
Also, the map $\theta:\tilde{G}\rightarrow V^\prime$ must be unramified, so $V^\prime$ admits the structure of an algebraic group in such a way that $\pi:V^\prime\rightarrow G$ is an isogeny.  
 \cvd

\medskip

To Prove Theorem 1.1, we begin by formulating (and proving)  a weaker version.
  \medskip

\noindent{\bf Proposition 5.3}. {\it Assume (a) and (b) of Theorem 1.1. Then there exists an algebraic group $G^\prime$, an isogeny $p: G^\prime\rightarrow G$  and a rational map $\theta:G^\prime\rightarrow X$, all defined over $\kappa$, such that for all $g^\prime$ in its domain,
$$
p(g^\prime)(\theta(g^\prime))=\theta(g^\prime).
$$
}
\smallskip

\noindent{\it Proof}.
We let $V\subset X\times G$ be the variety of fixed points for the given action of $G$ to $X$: 
$$
V:=\{(x,g)\in X\times G\, :\, g(x)=g\}.
$$
It is endowed with a projection $\pi:V\rightarrow G$. The hypothesis $(a)$ assures that $\pi$ is dominant, and hypothesis $(b)$ assures that its generic fiber is finite (say of degree $n$). The hypothesis $(a)$ of Theorem 1.1 states that $\Gamma\subset\pi  (V(\kappa))$. Then Proposition 5.1 provides the existence of an unramified covering $p: G^\prime\rightarrow G$ and a morphism $\theta^\prime:G^\prime\rightarrow V$ such that $\pi \circ\theta=p$. Letting $\pi_1: V\rightarrow X$ be the projection on the first factor, we obtain the rational map $\theta:=\pi_1 \circ\theta^\prime$, which has the required property.
\cvd
\medskip

Hence the first conclusion of Theorem 1.1 would be proved once we show that $G^\prime$ can be taken to be $G$ and $p$ to be the identity map. This is the content of the next result, of purely geometric nature.
\medskip

\noindent {\bf Proposition 5.4}. {\it Let $G$ be an algebraic group, $X$ an algebraic variety, both defined over a field $\kappa$ of characteristic zero. Suppose $G$ acts $\kappa$-morphically over $X$ in such a way that some element of $G$ has only finitely many fixed points in $X(\bar{\kappa})$. Suppose there exists an algebraic group $G^\prime$, an isogeny $p:G^\prime\rightarrow G$ and a rational map $\theta:G^\prime\rightarrow X$, defined over $\kappa$, such that $p(g)(\theta(g))=\theta(g)$ for all $g$ in the domain of $\theta$. Then such a map $\theta$ is constant on the fibers of $p$, i.e. is of the form $\theta=\omega\circ p$ for a rational map $\omega :G\rightarrow X$. Automatically $\omega$ satisfies $g(\omega(g))=\omega(g)$, for all $g$ in its domain.
}
\medskip

The proof of Proposition 5.4 will make use of analytic methods. Recall that the Lie algebra $\g$ of $G^\prime$ (which is isomorphic to the Lie algebra of $G$) is endowed with the exponential map $\exp:\g\rightarrow G^\prime(\C)$, whose image is an open neighborhod of the origin (possibly the whole Lie group $G^\prime(\C)$). An intermediate result toward the proof of Proposition 5.4 is
\medskip

\noindent{\bf Lemma 5.5}. {\it Under the hypotheses of Proposition 5.7, there exists a set $U \subset G^\prime(\C)$, open in the euclidean topology and dense in the Zariski topology,  such that for each $g^\prime \in U$, there exists an element $\alpha\in\g$ such that $g^\prime=\exp(\alpha)$ and the map $\C\ni t\mapsto\theta(\exp(t\alpha))$ is constant. 
}
\smallskip

\noindent{\it Proof}. Let $n$ be the generic number of fixed points for an element $g\in G$ in $X(\C)$ (such a number is finite by assumption). Let $U\subset G^\prime(\C)$ be the set of elements $g^\prime\in G^\prime$ with the properties that:
\smallskip

\item {(1)} there exists a point $\alpha$ in the Lie algebra of $G^\prime$ such that $g^\prime=\exp(\alpha)$;

\item {(2)}  $\exp(t\alpha)$  belongs to the domain of $\theta$ for all but finitely many $t\in\C$;   

\item {(3)} each element $g$ of $G$ of the form $g=p(\exp(t\alpha))$, for all but finitely many $t\in\C$, fixes at most $n$ points in $X(\C)$;
\smallskip

\noindent clearly $U$ is Zariski-dense and open in the euclidean topology. 

Let $g^\prime\in U$ and let $\alpha\in\g$ with $\exp(\alpha)=g^\prime$. Let us show that  the map $t\mapsto\exp(t\alpha)$ is constant, as claimed in the Lemma.
Assume by contradiction it is not so. Then, by continuity,  the  set   $\exp(t\alpha)$ for $t\in\Q$  would be an infinite set. Now, take $n+1$ rational points $t_1,\ldots,t_{n+1}\in\Q$  such that their images  $\theta(\exp(t_i\alpha))$ are pairwise distinct. There exist integers $d_1,\ldots,d_{n+1}$ such that the products $ d_it_i$ are all equal and none of them belong to the exceptional finite set of complex numbers $t$ such that $p(\exp(t\alpha))$ has more then $n$ fixed points. Neverthless, the point  $p(\exp(d_it_i\alpha))$, which is independent of $i$ and is a power of $p(\exp(t_i\alpha))$ for all $i=1,\ldots,n+1$, leaves fixed all the fixed points for $p(\exp(t_i\alpha))$, for $i=1,\ldots,n+1$. This contradiction proves the lemma.
\cvd

\medskip

\noindent{\it Proof of Proposition 5.4}. 
Let $Z\subset G^\prime$ be the kernel of the projection $p:G^\prime\rightarrow G$. It is a finite (normal) subgroup contained in the center of $G^\prime$; let $e$ be its order. Recall that $ \theta$ associates to a point $g^\prime\in G^\prime$ (in its domain) a fixed point in $X$ for the element $p(g^\prime)$ of $G$. We have to show that this fixed point is the same for $g^\prime$ and $g^\prime\cdot z$, whenever $z\in Z$ (at least for all $g^\prime$ in an open dense subset of $G^\prime$). It suffices to consider the points $g^\prime$ such that $g^\prime z\in U$ for all $z\in Z$ (they form a dense set in the Zariski topology, since $Z$ is finite). Let then   $g^\prime=\exp(\alpha)\in U$ be an element with this property. Then $\theta(g^\prime)=\theta({g^\prime}^e)$, since $\theta$ is constant on the one-parameter subgroup $t\mapsto\exp(t\alpha)$. For the same reason, we also have $\theta((g^\prime z)^e)=\theta(g^\prime z)$; on the other hand ${g^\prime}^d=(g^\prime z)^e$, since $z$ is central of order finite dividing $e$; so we have $\theta(g^\prime)=\theta(g^\prime z)$ as wanted.
\cvd

\medskip

\noindent{\it Proof of Theorem 1.1}. Conclusion $(i)$ of the Theorem follows immediately combining its the weak form given in Proposition 5.3 with  Proposition 5.4. So we know there is a $\kappa$-morphism $\omega:U\rightarrow X$, where $U\subset G$ is an open dense subset of $G$, satisfying (ii) of Theorem 1.1. Let us prove that, under the assumption that $X$ is projective, every element of $G(\kappa)$ has a rational fixed point. Let $g\in G(\kappa)$. We can find a smooth curve $\c\subset G$, defined over $\kappa$, passing through $g$ and not lying entirely in the complement of $U$. Then the restriction of $\omega$ to $\c$ can be continued to every point of $\c$, since it is a map from a smooth curve to a projective variety. The value of such continuation at the point $g$ is then a rational fixed point for $g$.\cvd
\medskip

{\it Proof of Theorem 1.2}. We shall apply Theorem 1.6, in its equivalent form given in Proposition 5.1. 
\medskip

\noindent Let $V\subset\Gm^r\times G$ be the variety defined by the condition
$$
(\lambda_1,\ldots,\lambda_r,g)\in V\qquad {\rm if\ and\ only\ if}\qquad (T-\lambda_1)\cdots(T-\lambda_r)|P(T,g),
$$
where as usual, for $g\in G\subset\GL_N$, $P(T,g)\in\kappa[T]$ denotes the characteristic polynomial of the matrix $g$. Denoting by $\pi:V\rightarrow G$ the projection onto the factor $G$, hypothesis $(i)$ of Theorem 1.2  garantees that $\Gamma\subset\pi(V(\kappa))$. Hence, by Proposition 5.1 there exist an algebraic group $\tilde{G}$, an isogeny $p:\tilde{G}\rightarrow G$ and rational functions $\tilde{\chi}_1,\ldots,\tilde{\chi}_r\in\kappa(\tilde{G})$ such that  for all $g^\prime$ in an open dense set of $\tilde{G}$, the polynomial $(T-\tilde{\chi}_1(g^\prime))\cdots(T-\tilde{\chi}_r(g^\prime))$ divides the characteristic polynomial of $\pi(g^\prime)$. To obtain the full Theorem 1.2, it remains to prove that: (1) the rational functions $\tilde{\chi}_1,\ldots,\tilde{\chi}_r$ are in fact regular; (2) one can take for $\tilde{G}$ the group $G$ and for $p:\tilde{G}\rightarrow G$ the identity; (3) the (regular) functions $\tilde{\chi}_1,\ldots,\tilde{\chi}_r$ are character of $G=\tilde{G}$.
\smallskip

\noindent  To prove our first claim, just observe that each function $\tilde{\chi}_i\in\kappa(\tilde{G})$ satisfies the monic equation $P(\tilde{\chi},p(g))=0$ over $\kappa[\tilde{G}]$, so is integral over the ring $\kappa[\tilde{G}]$. Since $\tilde{G}$ is smooth,  the corresponding regular function ring $\kappa[\tilde{G}]$ is integrally closed, so $\tilde{\chi}$ is a regular function. 

\noindent To prove both our second and third claims, we shall make use of Lemma 4.9. Since $\tilde{\chi}_1,\ldots,\tilde{\chi}_r$ are regular and never vanishing, they are characters of $\tilde{G}$, by Lemma 4.9. Let now $\tilde{\chi}\in\{\tilde{\chi}_1,\ldots,\tilde{\chi}_r\}$ be one of them.  We shall show that $\tilde{\chi}$ is constant in the pre-image $p^{-1}(g)$ of each point $g\in G$, thus proving that one can take $\tilde{G}=G$ and for $p$ the identity map, as wanted. Let $Z\subset\tilde{G}$ be the kernel of $p$; it is a finite central subgroup of order equal to the degree of $p$.  From $\tilde{\chi}(zg^\prime)=\tilde{\chi}(z)\tilde{\chi}(g^\prime)$, valid for all $z\in Z$, $g^\prime\in\tilde{G}$, it follows that the ratio $\tilde{\chi}(zg)/\tilde{\chi)}(g)$ is a fixed root of unity $\tilde{\chi}(z)$. Now, since $\tilde{\chi}(z)$ is an eigenvalue of $p(z)=1_G$, it follows that $\tilde{\chi}(z)=1$, so $\tilde{\chi}(zg)=\tilde{\chi}(g)$ as wanted.
\cvd
\medskip

\noindent{\it Proof of Corollary 1.3}. Let $G$ be the Zariski closure of the group $\Gamma$ in $\GL_N$, let $G^0$ be the neutral component of $G$ and put   $\Gamma^0:=\Gamma\cap G^0$. The subgroup $\Gamma^0\subset\Gamma$ has finite index in $\Gamma$ and is Zariski-dense in the connected algebraic group $G^0$; we shall prove it is solvable, obtaining the Corollary. The hypothesis of Theorem 1.2 are satisfied with $r=N$, with $G^0$ instead of $G$ and $\Gamma^0$ instead of $\Gamma$. Then by conclusion $(ii)$ of Theorem 1.1 there exists an algebraic-group homomorphism ${\bf \chi}:G^0\rightarrow\Gm^N$ whose kernel is composed of matrices all of whose eigenvalues are equal to $1$. Such a subgroup is known to be solvable; since the group $\Gm^N$ is also solvable, it follows that $G^0$ is solvable, hence so is $\Gamma^0$.\cvd

\medskip

\noindent{\it Proof of Corollary 1.4}. The equivalence between conditions $(i), (ii)$, and $(iii)$ is contained in Theorem 1.2. It then suffices  just to prove just that $(iii)$ implies $(iv)$ and $(iv)$ implies $(ii)$. 

Let us suppose that $(iii)$ holds, so there is a character $\chi:G\rightarrow \Gm$, defined over $\kappa$, such that $\chi(g)$ is for every $g\in G$ an eigenvalue of $g$. Let $r\in\{1,\ldots,N\}$ be the ``generic" value for the dimension of the kernel of $g-\chi(g)\cdot\uno$, where as usual $\uno$ denotes the unit-matrix in $\GL_N$. This number $r$ is also the minimal of such dimensions, for $g\in G(\bar{\kappa})$. 
Then we can define a  rational map $\bar{\omega}: G\rightarrow\F(r;N)$, to the Grassmanniann of $r$-dimensional subspaces in $\Ga^N$, by sending $g$ to $\bar{\omega}(g):={\rm ker} (g-\chi(g)\uno)$. Let $H$ be any subspace of codimension $r-1$ in $\Ga^N$, defined over $\kappa$, intersecting transversally at least one $r$-dimensional subspace of the form $\bar{\omega}(g)$.  Then the map $\omega: g\mapsto\bar{\omega}(g)\cap H$ sends a generic element of $G$ to a line in $\Ga^N$, i.e. to a point of $\Pr_{N-1}$; clearly such point is fixed for the projective automorphism induced by $g$. Hence we have proved $(iv)$ (assuming $(iii)$).

Let us now assume $(iv)$ and want to prove $(ii)$. It suffices to remark that a fixed point certainly exists for all $g$ in an open dense set $U$ of $G$ (where $\omega$ is well-defined). For any other point $g\in G(\kappa)$, letting $\c$ be a smooth curve on $G$ passing to $g$ and not lying entirely on $G\setminus U$, the restriction of the rational map $\omega$ to $\c$ can be continued to the whole curve $\c$, hence in particular to the point $g$; its value in $g$ provides a fixed point for $g$.
\cvd

\noindent{\it Proof of Theorem 1.5}. This is an easy consequence of Theorem 1.1. Actually,
the condition that some matrix of $G$ has distinct eigenvalues, appearing in the hypothesis of Theorem 1.5, assures that the natural action of $G$ on the flag variety $\F(r_1,\ldots,r_h;N)$ satisfies  condition $(b)$ in Theorem 1.1. Note that flag varieties are projective, so both conclusions of Theorem 1.1 hold. Now,  $(i)$  of Theorem 1.1 coincides exactly with $(ii)$ of Theorem 1.5, while conclusion $(ii)$ of Theorem 1.1 gives $(i)$ of Theorem 1.5. 
It only remains to prove the stronger conclusion in the particular case $(r_1,\ldots,r_h)=(1,\ldots,N-1)$, of the maximal flag variety. In this case we can apply Corollary 1.3, obtaining that $G$ is solvable. Then the Lie-Kolchin Theorem [Bo2, 10.5] assures the existence of a fixed complete flag for the whole group $G$.\cvd

\medskip

\noindent{\it Proof of Corollary 1.11}.  We follow closely [Se, \S 9.2], in particular its proof of Proposition 2 therehin. 
We first show that the points $\gamma\in G$ such that the Galois group (over $\kappa$) of the polynomial $P(T,\gamma)$ is {\it not} isomorphic to ${\cal G}$ is a $\kappa$-thin set. Let $V\subset\Gm^N\times G$ be the variety  
$$
V:=\{(\lambda_1,\ldots,\lambda_N,g)\in\Gm^N\times G \, :\, P(T,g)=(T-\lambda_1)\cdots(T-\lambda_N)\}.
$$ 
It is irreducible if and only if the characteristic polynomial of $G$ is irreducible.
Let $\pi:V\rightarrow G$ be the projection on the second factor. It is a (possibly disconnected) Galois cover of $G$ with Galois group isomorphic to ${\cal G}$. Let $V^{(1)},\ldots,V^{(r)}$, ($r\leq N$) be the irreducible components of $V$; each $V^{(j)}$ defines a connected cover of $G$, where the covering map is naturally the restriction of $\pi$ to $V^{(j)}$. Let ${\cal G}^{(j)}$ be the subgroup of ${\cal G}$ formed by the automorphisms acting trivially on $V^{(j)}$. The cover $\pi_{|V^{(j)}}: V^{(j)}\rightarrow G$  is  Galois with automorphism group ${\cal G}/{\cal G}^{(j)}$; in particular it has degree $>1$ whenever ${\cal G}^{(j)}\neq{\cal G}$.
Let now $\gamma\in G(\kappa)$ be a given matrix. Saying that the Galois group of the splitting field of $\gamma$ is not isomorphic to ${\cal G}$ amounts to the existence of an index $j\in\{1,\ldots,r\}$ such that: (1) ${\cal G}^{(j)}\neq {\cal G}$, (2) there exists a subgroup $H$ with ${\cal G}^{(j)}\subset H\subset {\cal G}$, $H\neq {\cal G}$, such that $H$ fixes a Galois invariant subset of the fiber of $\gamma$ in $V^{(j)}$. Consider now, for each index $j=1,\ldots,r$, the (possibly empty) set   $\{H_i^{(j)}\}_{i\leq i(j)}$ of subgroups ${\cal G}^{(j)}\subset H_i^{(j)}\subset{\cal G}$ with $H_i^{(j)}\neq{\cal G}$; there are no such subgroups if ${\cal G}^{(j)}={\cal G}$, i.e. if $\pi: V^{(j)}\rightarrow G$ is an isomorphism. 
 
Let $W$ be the union of the varieties of the form $V^{(j)}/H_i^{(j)}$. Consider the induced projection $\pi:W\rightarrow G$. By construction, it has no sections. Then the set of $\gamma\in G(\kappa)$ whose splitting field has a Galois group {\it not} isomorphic to ${\cal G}$ is just the image $\pi(W(\kappa))$ of the rational points in $W$, hence a $\kappa$-thin set. By Proposition 5.3, if such a set contains a Zariski-dense subgroup, then there exists a covering $p: \tilde{G}\rightarrow G$ and a rational map $\theta:\tilde{G}\rightarrow W$ with $p=\pi_W\circ\theta$. As in the proof of Theorems 1.1 and 1.2, it is easy to see that one can then choose $G^\prime=G$ and $p$ the identity map, concluding that $\pi$ admits a section, which we excluded.
\cvd  
\bigskip

\noindent {\bf References}. 
\medskip

\item{[Be]} {\teo J. Bernik}, On groups and semigroups of matrices with spectra in a finitely generated field, {\it Linear and Multilinear Algebra} {\bf 53,4} (2005), 259-267.\smallskip

\item{[B-O]} {\teo J. Bernik, J. Okni\'nski}, On semigroups of matrices with eigenvalue $1$ in small dimension, {\it Linear Algebra and its Applications}, {\bf 405} (2005), 67-73.\smallskip

\item{[B-M-Z]} {\teo E. Bombieri, D. Masser, U. Zannier}, Intersecting a curve with algebraic subgroups of multiplicative groups, {\it International Math. Research Notices} {\bf 20} (1999), 1119-1140.
\smallskip

\item{[Bo1]} {\teo A. Borel}, {\it Introduction aux groupes arithm\'etiques}, Hermann, Paris (1969).\smallskip

\item{[Bo2]} {\teo A. Borel}, {\it Linear Algebraic Groups}, 2nd Edition, GTM 126, Springer Verlag, 1997. \smallskip

\item{[C-R]} {\teo C. W. Curtis, I. Reiner}, {\it Representation theory of finite groups and associative algebras}, John Wiley \& Sons, 1962.\smallskip

\item{[D]} {\teo P. D\`ebes},  On the irreducibility of the polynomial $P(t^m,Y)$, {\it J. Number Theory}, {\bf 42} (1992), 141-157.\smallskip

\item{[D-Z]} {\teo R. Dvornicich, U. Zannier}, Cyclotomic Diophantine Problems (Hilbert Irreducibi\-li\-ty and Invariant Sets for Polynomial Maps), preprint (2006). 
\smallskip

\item{[F-Z]} {\teo A. Ferretti, U. Zannier},  Equations in the Hadamard ring of rational functions, {\it Annali Scuola Normale Sup.} (to appear).\smallskip

\item{[H]} {\teo D. Hilbert}, Ueber die Irreducibilit\"at ganzer rationaler Functionen mit ganzzahligen Coefficienten, J. reine ang. Math. {\bf 110} (1892), 104-129.
\smallskip

\item{[L]} {\teo S. Lang}. {\it Fundamentals of Diophantine Geometry}, Springer Verlag, 1984.\smallskip

\item{[Lau]} {\teo M. Laurent}, Equations diophantiennes exponentielles et suites r\'ecurrentes lin\'eaires II, {\it Journal of Number Theory}, {\bf 31} (1988), 24-53. \smallskip

\item{[M]} {\teo D. W. Masser}, Specializations of finitely generated subgroups of abelian varieties, {\it Transactions Am. Math. Soc.}, {\bf 311} (1989), 413-424.\smallskip


\item{[vdP]} {\teo A. van der Poorten}, Some facts that should be better known, especially about rational functions, in Number Theory and Applications (Banff, AB 1988), 497-528, Kluwer Acad. Publ., Dordrecht, 1989. \smallskip

\item{[P-R 1]} {\teo G. Prasad, A. Rapinchuk}, Existence of irreducible $\R$-regular elements in Zariski-dense subgroups, {\it Math. Research Letters} {\bf 10} (2003), 21-32.\smallskip

\item{[P-R 2]} {\teo G. Prasad, A. Rapinchuk}, Zariski-dense subgroups and transcendental number theory, {\it Math. Research Letters} {\bf 12} (2005), 239-249.\smallskip

\item{[R]} {\teo R. Rumely}, Notes on van der Poorten proof of the Hadamard quotient theorem II, in: S\'eminaire de Th\'eorie des Nombres de Paris 1986-87, {\it Progress in Mathematics}, Birkh\"auser 1988.\smallskip

\item{[Se]} {\teo J.-P. Serre}, {\it Lectures on the Mordell-Weil Theorem}, 3rd Edition, Vieweg-Verlag, 1997. \smallskip

\item{[Sch]} {\teo W. M. Schmidt}, {\it Linear Recurrence Sequences and Polynomial-Exponential Equations}, in Amoroso \& Zannier (ed.) {\it Diophantine Approximation}, Proceedings of the C.I.M.E. Conference, Cetraro 2000, Springer LNM 1829 (2003).
\smallskip

\item{[Z1]} {\teo U. Zannier}, A proof of Pisot's $d$-th root conjecture, {\it Annals of Math.} {\bf 151} (2000), 375-383.
\smallskip

\item{[Z2]} {\teo U. Zannier}, {\it Some applications of diophantine approximations to diophantine equations}, Forum Editrice, Udine  2003. 
\bigskip

$${}$$

\bigskip

{\teo Dipartimento di Matematica e Informatica}

Via delle Scienze, 206 - 33100 Udine

{\tt corvaja@dimi.uniud.it}

\bye